\documentclass[a4paper,11pt]{article}
\usepackage{amssymb}
\usepackage{amsmath}
\usepackage{amsfonts}
\usepackage[mathscr]{eucal}
\usepackage{enumerate}
\usepackage{amscd}
\usepackage{indentfirst}
\usepackage{enumerate}
\usepackage{txfonts, textcomp}
\newtheorem{thm}{Theorem}[section]
\newtheorem{prop}[thm]{Proposition}
\newtheorem{dfn}[thm]{Definition}
\newtheorem{lem}[thm]{Lemma}

\newtheorem{ex}{Example}[section]
\newtheorem{rmk}{Remark}[section]
\numberwithin{equation}{section}
\allowdisplaybreaks[4]

\title{GEOMETRIC QUANTIZATION OF DIRAC MANIFOLDS}
\author{Yuji HIROTA}
\date{}
\begin{document}
\maketitle 

\begin{abstract}
We define prequantization for Dirac manifolds to generalize known procedures for Poisson and (pre) symplectic manifolds 
by using characteristic distributions obtained from 2-cocycles associated to Dirac structures. 
Given a Dirac manifold $(M,D)$, we 
construct a Poisson structure on the space of admissible functions on $(M,D)$ and a representation of 
the Poisson algebra to establish the prequantization condition of $(M,D)$ in terms of a Lie algebroid cohomology. 
Additional to this, we introduce a polarization 
for a Dirac manifold $M$ and discuss procedures for quantization in two cases where 
$M$ is compact and where $M$ is not compact.
\end{abstract}

\section{Introduction}

In classical mechanics, a state of a system is described by a pair of position and momentum, 
and the time evolution of the system is controlled by Hamilton's equation. On the other hand, in quantum 
mechanics, a quantum state is given as a point in a complex Hilbert space, and the time evolution of 
quantum system is described by Shr$\ddot{\rm o}$dinger's equation. In addition, physical quantities such as 
position, momentum and energy are given by self-adjoint operators on the Hilbert space describing the quantum 
system, and the time evolution of the physical quantity is defined by Heisenberg's equation. 
For instance, we consider the motion of $N$-particles in $\mathbb{R}^3$ and let $(q_j,p_j)~(j=1,2,3)$ be the 
coordinates of position and momentum variables. Then, the space $C^\infty(\mathbb{R}^{6N})$ of smooth functions on 
the phase space $\mathbb{R}^{6N}\simeq \mathbb{R}^{3N}\times\mathbb{R}^{3N}$ is Poisson algebra by 
\[
\{F,G\}:=\sum_j\Bigl(\,\frac{\partial G}{\partial p_j}\frac{\partial F}{\partial q_j}-
 \frac{\partial F}{\partial p_j}\frac{\partial G}{\partial q_j}\,\Bigr).
\] 
This leads to the following relations: 
\begin{equation*}\label{eqn}
 \{q_j,\,q_k\} = \{p_j,\,p_k\} = 0,\qquad \{q_j,\,p_k\} = \delta_{j,k}\quad (\forall j,k=1,2,\cdots,3N). 
\end{equation*}
Defining self-adjoint operators $\hat{q}_j$ and $\hat{p}_j$ in a Hilbert space $L^2(\mathbb{R}^{3N})$ as 
$\hat{q}_j:= {q_j}_{\cdot},\,\hat{p}_j:=-\sqrt{-1}\hbar\,\frac{\partial}{\partial q_j}$ for 
positions $q_j$ and momenta $p_j$, one can 
get the following relations, called the canonical commutation relations: 
\begin{equation*}\label{eqn2}
 [\hat{q}_j,\,\hat{q}_k] = [\hat{p}_j,\,\hat{p}_k] = 0,\qquad 
 [\hat{q}_j,\,\hat{p}_k] = \mathrm{\bf i}\hbar\,\delta_{j,k}\quad (\forall j,k=1,2,\cdots,3N), 
\end{equation*}
where $[\cdot,\cdot]$ means a commutator $[A, B]:=AB-BA$ for operator algebras $A, B$.  
In other words, one can obtain quantum objects from classical objects 
by choosing a Hilbert space proper for corresponding to the classical theory and 
by constructing self-adjoint operators on it. 

Such a procedure to determine a quantum theory which corresponds to a given classical theory 
is called a quantization. Mathematically, a quantization is a procedure to construct a representation of a 
Poisson algebra on the space of functions on a proper Hilbert space for a given manifold. 
Quantization is a meaningful subject for the study of mathematics and physics, and is an extremely 
interesting one at which mathematics and physics intersect as well. 
It is known that there are several kinds of quantizations, such as 
canonical quantization, Feynman's path integral quantization, geometric quantization, 
Moyal quantization, Weyl-Wigner quantization and so on. 
Among those, geometric quantization consists of two procedures: prequantization and polarization. 
Prequantization assigns to a given symplectic manifold $S$ a Hermitian line bundle $L\to S$ with a 
connection whose curvature 2-form is the symplectic structure. 
Then, a Poisson subalgebra of $C^\infty(S)$ acts faithfully on the space 
$\varGamma^\infty(S,L)$ of smooth sections of $L$. 
On the other hand, a polarization is the procedure which reduces $\varGamma^\infty(S,L)$ to 
a subspace $\mathcal{A}\subset \varGamma^\infty(S,L)$ appropriate for physics 
so that a subalgebra of $C^\infty(S)$ may still act on $\mathcal{A}$. 

The study of geometric quantization in symplectic geometry goes back to the theory by 
B. Kostant and J. Souriau (see \cite{Kqua70, Sstr69}). Later, a target for quantization was extended to 
a presymplectic manifold, and its quantization was studied by many researchers 
in \cite{GSont81, Gpre80, Vgeo83}.  
After that, J. Huebschmann extended the target from a (pre)symplectic manifold to a Poisson manifold and 
investigated algebraically its quantization in \cite{Hpoi90}. Besides, I. Vaisman 
studied the quantization of a Poisson manifold in terms of Hermitian line bundles in \cite{Vont91}, 
and D. Chinea, J. Marrero and M. de Leon did it in terms of $S^1$-bundles in \cite{CMLpre96}. 
In the case where the target is a twisted Poisson manifold, its geometric quantization 
is studied by F. Petalidou in \cite{Pont07}. Lastly, in \cite{WZvar05}, A. Weinstein and M. Zambon studied 
a prequantization of Dirac manifolds which is a generalization of (pre) symplectic and Poisson manifolds 
in terms of Dirac-Jacobi structures appeared in \cite{Wcon00, IMlie02}. 

Dirac manifolds were introduced by T. Courant for the purpose of unifying approaches to the 
geometry of Hamiltonian vector fields and their Poisson algebras, 
which are thought of generalizations of both presymplectic manifolds and Poisson manifolds (see \cite{Cdir90}). 
The purpose of this paper is to unify known prequantization approach for Poisson and (pre) symplectic manifolds 
by introducing a prequantization procedure for Dirac manifolds. 
Using characteristic distributions from a 2-cocycle associated to $(M,D)$, 
we define a Poisson bracket for admissible functions on $(M,D)$ and 
give a representation of the Poisson algebra in terms of a connection theory of Lie algebroids 
to describe a condition for prequantization of $(M,\,D)$. 
Our approach to prequantization for Dirac manifolds is different from the one discussed in \cite{WZvar05}. 

The paper is organized as follows: In Section 2, we review the fundamentals of Dirac manifolds. 
In Section 3, after reviewing the Lie algebroid cohomology and the connection theory of 
Lie algebroids, we introduce the first Dirac-Chern class of line bundles over Dirac manifolds 
and show that it does not depend on a choice of connections. 
In Section 4, we introduce a prequantization for Dirac manifolds. 
We define a Poisson structure on the space of admissible functions 
associated to a singular distribution for a given Dirac manifold $(M,D)$ and 
construct a map from the Poisson algebra to the space of sections of a complex line bundle over $(M,D)$. 
We provide the necessary and sufficient condition for the map to be a representation of the Poisson algebra. 
Lastly, we formulate the condition for prequantization of $(M,D)$ to be realized in terms of Lie algebroid 
cohomology for $(M,D)$. In Section 5, we introduce polarizations for Dirac manifolds 
and develop the quantization process of them, basing on the discussion in Section 4. 
\vspace{0.3cm}

Throughout the paper, every smooth manifold is assumed to be paracompact, and all maps are assumed to be smooth. 
We denote by $\varGamma^\infty(M,E)$ the space of 
smooth sections of a smooth vector bundle $E\to M$. 
Especially, if $E=TM$, 
we often write $\mathfrak{X}\,(M)$ for $\varGamma^\infty(M,TM)$. 
We use $\Omega^k(M)$ and $\mathfrak{X}^k(M)$ for $\varGamma^\infty(M, \wedge^kT^*M)$ and 
$\varGamma^\infty(M, \wedge^kTM)$, respectively. 
\section{Dirac manifolds}

\subsection{Definition and examples}
Let $M$ be a finite dimensional smooth manifold. We define symmetric and skew-symmetric 
operations on the vector bundle $\mathbb{T}M:=TM\oplus T^*M$ over $M$ as 
\[
\langle\,X\oplus \xi,\,Y\oplus\eta\,\rangle_+ := \frac{1}{2}\,\bigl\{\xi\,(Y) + \eta\,(X)\bigr\} 
 \in C^\infty(M)
\]
and 
\[
\llbracket X\oplus\xi,\,Y\oplus\eta\rrbracket 
:= [X,Y]\oplus (\mathcal{L}_X\eta - \mathrm{i}_Yd\xi) \in \varGamma^\infty(M, \mathbb{T}M)
\]
for all $X\oplus\xi, Y\oplus\eta \in \varGamma^\infty(M, \mathbb{T}M)$. 
Here $\mathcal{L}_X\eta$ stands for the Lie derivative of $\eta$ 
by $X$ and $\mathrm{i}_Yd\xi$ for the contraction of $d\xi$ with $Y$. 
A subbundle $D\subset \mathbb{T}M$ is called a Dirac structure on $M$ 
if the following conditions are satisfied: 
\begin{enumerate}[\quad(D1)]
 \item $\left.\langle\cdot,\,\cdot\rangle_+\right|_D \equiv 0$; 
 \item $D$ has rank equal to $\dim M$;
 \item 
  $\llbracket\varGamma^\infty(M, D),\,\varGamma^\infty(M, D)\rrbracket\subset 
  \varGamma^\infty(M, D)$.
\end{enumerate}
A smooth manifold $M$ together with 
Dirac structure $D\subset \mathbb{T}M$ is called a Dirac manifold, denoted by $(M, D)$. 
In addition to the natural pairing $\langle\cdot,\,\cdot\rangle_+$, one defines a 
skew-symmetric pairing $\langle\cdot,\,\cdot\rangle_-$ as 
\[
\langle\,X\oplus\xi,\,Y\oplus\eta\,\rangle_- := \frac{1}{2}\,\{\xi\,(Y) - \eta\,(X)\}  \in C^\infty(M). 
\]
The following formula can be shown by direct calculation and Condition (D3). 
\begin{lem}\label{sec2:integrability}
It holds that 
\begin{equation}\label{sec2:D3}
 (\mathcal{L}_{X_1}\xi_2)(X_3) \;+\; (\mathcal{L}_{X_2}\xi_3)(X_1) \;+\; 
(\mathcal{L}_{X_3}\xi_1)(X_2) \;=\; 0 
\end{equation}
for any $X_1\oplus \xi_1,\,X_2\oplus \xi_2,\, X_3\oplus \xi_3\in\varGamma^\infty(M,D)$. 
\end{lem}

Given a Dirac structure $D\subset \mathbb{T}M$, there are natural 
projections 
\[
(\rho_{TM})_p:=\mathrm{pr}_1|_{D_p}: D_p\to T_pM\qquad \text{and}\qquad
(\rho_{T^*M})_p:=\mathrm{pr}_2|_{D_p}: D_p\to T_p^*M 
\]
for each point $p\in M$. 
The singular distribution $M\ni p\mapsto (\rho_{TM})_p(D_p)\subset T_pM$ is called the characteristic 
distribution. 
It is known that the characteristic distribution is integrable in the sense of \cite{DZpoi05} and \cite{Sorb73}. 
The corresponding singular foliation is called the characteristic foliation. 
For a discussion of singular distributions and the integrability, we refer to \cite{Sacc74a, Sacc74b} and \cite{Sorb73}. 
It is easy to check that 
\begin{equation}\label{sec2:kernel}
\ker\rho_{TM} = D\cap T^*M\qquad \text{and}\qquad
\ker\rho_{T^*M}= D\cap TM, 
\end{equation}
where $D\cap T^*M:=D\cap (\{\boldsymbol{0}\}\oplus T^*M),\,
D\cap TM:=D\cap (TM\oplus\{\boldsymbol{0}\})$. 
Here, we remark that, for each $p\in M$, $D_p\cap T_pM$ (resp. $D_p\cap T_p^*M$) are thought of as a subspace of 
either $T_pM\oplus T_p^*M$ or $T_pM$ (resp. $T_p^*M$). 
The following proposition is easily checked. 

\begin{prop}\label{sec2:prop}$\mathrm{(\cite{Cdir90})}$~
Given a Dirac manifold $(M, D)$, one has the characteristic equations 
\[
\rho_{TM}(D)= (D\cap T^*M)^\circ\qquad \text{\rm and}\qquad
\rho_{T^*M}(D)= (D\cap TM)^\circ, 
\]
where the symbol ${}^\circ$ stands for the annihilator. 
\end{prop}

For each $p\in M$, we define a bilinear map $\Omega_p$ on the subspace 
$(\rho_{TM})_p(D_p)\subset T_pM$ as 
\begin{equation}\label{sec2:2-form}
\Omega_p(X_p,\,Y_p) := \xi_p\,(Y_p)\quad (\forall Y_p\in(\rho_{TM})_p(D_p)),
\end{equation}
where $\xi_p$ is an element in $T_p^*M$ such that $X_p\oplus \xi_p\in D_p$. 
It is shown that $\Omega$ is well-defined, and is a presymplectic form on $\rho_{TM}(D)$ by Lemma \ref{sec2:integrability} 
and Proposition \ref{sec2:prop}~(see \cite{Cdir90}). 
The symbol $\Omega^\flat$ denotes 
the bundle map induced from $\Omega$. That is, $\Omega^\flat$ is the map 
$\Omega^\flat:\rho_{TM}(D)\to \rho_{TM}(D)^*$ 
which assigns $\Omega^\flat(X) = \Omega(X,\cdot)$ to $X\in \rho_{TM}(D)$. 
One easily finds that $\ker\Omega^\flat=D\cap TM$. 

In the same way, one also obtains a skew-symmetric tensor fields 
$\Pi:\rho_{T^*M}(D) \times \rho_{T^*M}(D) \to C^\infty(M)$ by 
\begin{equation}\label{sec2:bivector}
 \Pi_p(\xi_p,\eta_p) := \xi_p\,(Y_p)\quad (p\in M), 
\end{equation}
where $Y_p$ is a vector in $T_pM$ such that $Y_p\oplus \eta_p\in D_p$. 
The form $\Pi$ defines a map, denoted by $\Pi^\sharp$, from the subspace 
$\rho_{T^*M}(D)=(D\cap TM)^\circ$ to $(\rho_{T^*M}(D))^*$ by 
\[
 (\rho_{T^*M})_p(D_p)\ni\eta_p\longmapsto 
 \bigl\{\,\xi_p\mapsto \,\xi_p\,\bigl(\Pi_p^\sharp(\eta_p)\bigr):=\xi_p\,(Y_p)\,\bigr\}
 \in (\rho_{T^*M}(D))^* 
\]
for each $p\in M$. 
Letting $X\oplus \xi, Y\oplus \eta$ be smooth sections of $D$, we have 
$X = \Pi^\sharp(\xi), Y=\Pi^\sharp(\eta)$ and 
\[ 
\llbracket \Pi^\sharp(\xi)\oplus \xi,\,\Pi^\sharp(\eta)\oplus \eta\rrbracket
 = [\Pi^\sharp(\xi), \Pi^\sharp(\eta)]\oplus\{ \xi, \eta \} 
 \in \varGamma^\infty(M,D), 
\]
where $\{\xi, \eta\} := \mathcal{L}_{\Pi^\sharp(\xi)}\eta - \mathrm{i}_{\Pi^\sharp(\eta)}d\xi$. 
This implies that 
\[
 \Pi^\sharp(\{ \xi, \eta \}) = [\Pi^\sharp(\xi), \Pi^\sharp(\eta)]\quad 
 (\,\forall X\oplus \xi,\, Y\oplus \eta \in \varGamma^\infty(M,D)\,). 
\]

\begin{ex}\label{sec2:example}
Suppose that $M$ be a (pre)symplectic manifold with a (pre)symplectic form $\omega$. 
The 2-form $\omega$ induces the bundle map 
\[
\omega^\flat : \mathfrak{X}\,(M)\longrightarrow \Omega^1(M),\qquad 
X \longmapsto \mathrm{i}_X\omega . 
\]
One can obtain the subbundle $\mathrm{graph}\,(\omega^\flat)$ in $\mathbb{T}M$ as 
\[ 
\mathrm{graph}\,(\omega^\flat)_p := \{\,X_p\oplus \mathrm{i}_{X_p}\omega_p\in T_pM\oplus T_p^*M
\,|\,X_p\in T_pM\,\} \quad (p\in M)
\]
and can verify that $\mathrm{graph}(\omega^\flat)$ satisfies the three conditions 
{\rm (D1) -- (D3)} in the above. Therefore, $(M,\,\mathrm{graph}\,(\omega^\flat))$ is a Dirac manifold. 
Similarly, any symplectic manifold $M$ defines a Dirac structure on $M$. 
\end{ex}

\begin{ex}\label{sec2:example2}
Similarly to Example \ref{sec2:example}, any Poisson manifold $(P,\,\pi)$ defines a 
Dirac structure. Indeed, the 2-vector field $\pi$ induces the bundle map 
\[
\pi^\sharp : \Omega^1(P)\longrightarrow \mathfrak{X}\,(P),\qquad 
\alpha \longmapsto \{\,\beta\mapsto \pi\,(\beta, \alpha)\,\}. 
\]
and the subbundle $\mathrm{graph}\,(\pi^\sharp)$ given by 
\[
\mathrm{graph}\,(\pi^\sharp)_p := \{\,\pi^\sharp(\xi_p)\oplus \xi_p\in T_pP\oplus T_p^*P
\,|\,\xi_p\in T_p^*P\,\} \quad (p\in P). 
\] 
It can be easily verified that $(P,\,\mathrm{graph}\,(\pi^\sharp))$ is a Dirac manifold. 
\end{ex}

\begin{ex}\label{sec2:regular foliation}
 We let $F\subset TM$ be a regular distribution and assume 
 that $F$ is involutive. Then, the vector bundle $F\oplus F^\circ\to M$ is a Dirac structure on $M$, 
 where $F^\circ$ denotes the annihilator of $F$ in $T^*M$. 
\end{ex}

\begin{ex}
 Let $Q$ be a submanifold of a Dirac manifold $(M,\,D)$. If either $D_q\cap (T_qQ\oplus T_q^*M)$ or $D_q\cap (T_qQ)^\circ$
 has constant dimension at each point $q\in Q$, $Q$ has a Dirac structure $D_Q$ defined by 
 \[
   (D_Q)_q = \frac{D_q\cap (T_qQ\oplus T_q^*M)}{D_q\cap (T_qQ)^\circ}. 
 \]
\end{ex}

\begin{ex}\label{sec2:example,1-form}
Let $\eta$ be a 1-form on $M$. A subbundle 
\[
 (D_\eta)_p := \{\,X_p\oplus \mathrm{i}_{X_p}(d\eta)_p\in T_pM\oplus T_p^*M\,|\,X_p\in T_pM\,\} \quad (p\in M)
\]
satisfies the three conditions {\rm (D1) -- (D3)}. Therefore, $(M,\,D_\eta)$ is a Dirac manifold. 
\end{ex}

\begin{ex}\label{sec2:example, canonical}
 We define a 2-form on $T^*M\times\mathbb{R}$ as the pullback of the canonical symplectic form $\omega_0$ on 
 $T^*M$ by the projection $\mathrm{pr}_1$ on the first factor. Then, a subbundle
 \[
  (D_{\omega_0})_{(z,t)} = \Bigl\{\,\Bigl(X_z,\,f(z,t)\left.\frac{d}{dt}\right|_t\Bigr) \oplus \mathrm{pr}_1^*(\mathrm{i}_{X_z}\omega_0)
\,\bigm|\,\Bigl(X_z,\,f(z,t)\left.\frac{d}{dt}\right|_t\Bigr) 
   \in T_{(z,t)}(T^*M\times\mathbb{R})\,\Bigr\}\quad ((z,t)\in T^*M\times \mathbb{R})
 \]
 turns out to be a Dirac structure over $T^*M\times\mathbb{R}$ by noting that the vector field $X$ on $T^*M$ is $\mathrm{pr}_1$-related 
 to $\bigl(X,\,f\frac{d}{dt}\bigr)$ on $T^*M\times\mathbb{R}$. 
\end{ex}

\begin{ex}\label{sec2:example, almost cosymplectic}
 Let $(M,\,\omega,\,\eta)$ be an almost cosymplectic manifold. That is, $M$ is a $2k+1$-dimensional manifold equipped with 
 a 2-form $\omega$ and a 1-form $\eta$ such that $\eta\wedge\omega^k$ is a volume form on $M$. If $\omega$ is 
 closed, then a subbundle 
 \[
  (D_{\omega,\eta})_p = \{\,X_p\oplus \mathrm{i}_{X_p}(\omega + d\eta)\,|\,X_p\in T_pM\,\}\quad (p\in M)
 \]
 is a Dirac structure over $M$. 
\end{ex}

\begin{ex}
 A Jacobi manifold $(M,\,\pi,\, E)$ is a smooth manifold equipped with a bivector field $\pi$ and a vector field $E$ on $M$ 
 which satisfy $[\pi,\,\pi]_{\rm SN}=2E\wedge \pi$ and $[E,\,\pi]_{\rm SN}=0$, where $[\cdot,\,\cdot]_{\rm SN}$ denotes the 
 Schouten-Nijenhuis bracket. Let us consider a subbundle $L$ of $(TM\times\mathbb{R})\oplus(T^*M\times\mathbb{R})$ over $M$ by 
 \[
  L_p = \bigl\{\,(\pi^\sharp\alpha_p + f(p)E_p,\,-\alpha_p(E_p))\oplus (\alpha_p,\,f(p))\,
             |\,(\alpha_p,\,f(p))\in T_p^*M\times \mathbb{R}\,\bigr\}\quad (p\in M). 
 \]
 The subbundle $L$ gives rise to a Dirac structure $\tilde{L}\subset\mathbb{T}(M\times\mathbb{R})$ over $M\times \mathbb{R}$ by 
 \[
   \tilde{L}_{(p,t)} = \Bigl\{\,\bigl(\pi^\sharp\alpha_p + f(p)E_p,\,-\alpha_p(E_p)\partial_t\bigr)
  \oplus e^t\bigl(\alpha_p,\,f(p)(dt)_t\bigr)\,
  \bigm|\,(\alpha_p,\,f(p)\partial_t)\in T_p^*M\times T_t\mathbb{R} \,\Bigr\}\quad ((p,t)\in M\times \mathbb{R}),
 \]
 where $\partial_t = \left.\frac{\partial}{\partial t}\right|_t$~(see Section 5 in \cite{IMlie02}). 
\end{ex}

\begin{ex}
 Let $(M,\,\pi,\,E)$ be a Jacobi manifold of dimension $n$ and $z$ a point in $M$ where $E_z\ne \boldsymbol{0}$. 
 Suppose that $u_1,\cdots,u_{n-1}$ are functions on a neighborhood $U_z$ such that 
 $(du_1)_p,\cdots,(du_{n-1})_p$ are linearly independent at each $p\in U_z$ and $du_i\in E^\circ~(i=1,\cdots,n-1)$. 
 A subbundle 
 \[
  D_{\pi,E} = {\rm span}\,
     \bigl\{\,(\pi^\sharp(du_1)+u_1E)\oplus du_1,\cdots, (\pi^\sharp(du_{n-1})+u_{n-1}E)\oplus du_{n-1},\,E\oplus\boldsymbol{0}\,\bigr\}
 \]
 of $\mathbb{T}U_z$ is a Dirac structure over $U_z$~(see the subsection 4.2 in \cite{Cdir90}).
\end{ex}

\subsection{Admissible functions}

A smooth function $f$ on a Dirac manifold $(M, D)$ is said to be admissible 
if there exists a vector field $X_f\in \mathfrak{X}\,(M)$ such that 
$X_f\oplus df$ is a smooth section of $D$ (see \cite{Cdir90}). 
We note that the vector field $X_f$ is not uniquely determined. 
As easily checked by the case of a Dirac manifold induced by a presymplectic structure 
(see Example \ref{sec2:example}), $X_f$ is uniquely defined up to elements of $\ker\Omega$. 

\begin{ex}\label{sec2:example presymplectic}
Consider the presymplectic structure $\omega = dx_1\wedge dx_2 + dx_1\wedge dx_4$ on $\mathbb{R}^4$ 
and a function $f(x_1,x_2,x_3,x_4) = {x_1}^2 + k\,(x_2 + x_4)$ $(k\in\mathbb{R})$. 
Then, vector fields written in the form 
\[X = k\frac{\partial}{\partial x_1} + \varphi_1(\boldsymbol{x})\frac{\partial}{\partial x_2} 
 + \varphi_2(\boldsymbol{x})\frac{\partial}{\partial x_3} 
- (2x_1 + \varphi_1(\boldsymbol{x}))\frac{\partial}{\partial x_4}\qquad 
 \bigl(\,\varphi_1,\,\varphi_2\in C^\infty(\mathbb{R}^4)\,\bigr), 
\]
turn out to satisfy that $\omega^\flat(X)=df$. 
Therefore, $f$ is an admissible function on $(\mathbb{R}^4, \mathrm{graph}\,(\omega^\flat))$. 
\end{ex}

Given a Dirac manifold $(M,\,D)$, we denote the space of admissible functions on $(M,\, D)$ by 
$C^\infty_{\rm adm}(M,D)$. 
For any admissible functions $f,\,g\in C^\infty_{\rm adm}(M,D)$, one defines their bracket 
$\{f,\,g\}'$ as 
\begin{equation}\label{sec2:bracket}
\{f,\,g\}' := X_gf. 
\end{equation}
It can be shown that the bracket (\ref{sec2:bracket}) is both well-defined 
and skew-symmetric in the same way as
the case of $\Omega$. If $f,\,g$ are admissible, there exist vector fields 
$X_f$ and $X_g$ on $M$ such that $(X_f,\,df), (X_g,\,dg)\in\varGamma^\infty(M,\,D)$. 
Then, the simple computation yields 
\begin{equation*}
 \llbracket\, X_g\oplus dg, X_f\oplus df\,\rrbracket 
= (-[X_f,\,X_g])\oplus d\{f, g\}'\in\varGamma^\infty(M,\,D). 
\end{equation*}
This implies that the bracket $\{f,\,g\}'$, also, is admissible and satisfies the equation 
\begin{equation}\label{sec2:bracket2}
 X_{\{f,g\}'} \,+\, [X_f,\,X_g] \,=\, \boldsymbol{0}. 
\end{equation}
The next proposition can be shown by using (\ref{sec2:bracket2}) (see \cite{Cdir90}). 

\begin{prop}
 $(\,C^\infty_{\rm adm}(M,D),\,\{\cdot,\cdot\}'\,)$ forms a Poisson algebra. 
\end{prop}

\section{Lie algebroids}

\subsection{Basic terminology}

To carry out the procedure of prequantization for Dirac manifolds, 
the notion of cohomology for Dirac manifold is needed. 
Before proceeding the discussion, let us recall the definition of Lie algebroid and 
its cohomology. 
\begin{dfn}
A Lie algebroid over $M$ is a smooth vector bundle $A\to M$ 
with a bundle map $\sharp:A\to TM$, called the anchor map, and 
a Lie bracket $[\cdot,\,\cdot]$ on the space $\varGamma^\infty(M, A)$ of 
smooth sections of $A$ such that 
\begin{equation}\label{sec2:eqn1}
 [\alpha,\,f\beta] 
            = \bigl((\sharp\alpha) f\bigr)\,\beta \,+\, f[\alpha,\,\beta]
\end{equation}
for any $f\in C^\infty(M)$ and $\alpha,\,\beta\in \varGamma^\infty(M, A)$. 
\end{dfn}

A simple example is a tangent bundle $TM$ over a smooth manifold $M$: 
the anchor map $\sharp$ is the identity map, and the bracket $[\cdot,\cdot]$ 
is the usual Lie bracket of vector fields. This is called the tangent algebroid of $M$. 
As is well-known, Poisson manifolds define the structure of Lie algebroid 
on their cotangent bundles. 

\begin{ex}\label{sec3:cotangent algebroid}
(Cotangent algebroids)~
If $(P,\,\pi)$ is a Poisson manifold, then a cotangent bundle $T^*P$ is a Lie algebroid$\mathrm{:}$ 
the anchor map is the map $\pi^\sharp$ induced from $\varpi$, 
\[\pi^\sharp: T^*P\longrightarrow TP,\quad \alpha\longmapsto 
\bigl\{\, \beta\mapsto \langle\beta,\,\varpi^\sharp(\alpha)\rangle 
= \pi(\beta,\,\alpha)\, \bigr\}\]
and the Lie bracket is given by 
\begin{align*}
\{\alpha,\,\beta\} &:= \mathcal{L}_{\pi^\sharp(\alpha)}\beta
 - \mathcal{L}_{\pi^\sharp(\beta)}\alpha + d\bigl(\pi(\alpha,\,\beta)\bigr)\\
 &= \mathcal{L}_{\pi^\sharp(\alpha)}\beta - \mathrm{i}_{\pi^\sharp(\beta)}d\alpha 
\end{align*} 
as in the part immediately before the subsection 2.2. 
The Lie algebroid $(T^*P\to P,\,\{\cdot,\,\cdot\},\,\varpi^\sharp)$ is called a 
cotangent algebroid. 
\end{ex}

For other examples and the fundamental properties of Lie algebroids, 
see \cite{CWgeo99} and \cite{DZpoi05}. 

Let $(A_1\to M_1,\,[\cdot,\,\cdot]_1,\,\sharp_1)$ and 
$(A_2\to M_2,\,[\cdot,\,\cdot]_2,\,\sharp_2)$ be Lie algebroids. 
A Lie algebroid morphism from $A_1$ to $A_2$ is a vector bundle morphism 
$\Phi:A_1\to A_2$ with a base map $\varphi : M_1\to M_2$ which satisfies 
\begin{equation*}
   \sharp_2\bigl(\Phi(\alpha)\bigr) = \varphi_*\bigl(\sharp_1(\alpha)\bigr),\quad 
  \bigl(\forall \alpha\in \varGamma^\infty(M_1, A_1)\bigr),
\end{equation*}
and, for any smooth sections $\alpha,\,\beta\in \varGamma^\infty(M_1, A_1)$ written in the forms 
\begin{equation*}
\Phi\circ\alpha = \sum_i\xi_i\,(\alpha'_i\circ \varphi),\quad 
\Phi\circ\beta = \sum_j\eta_j\,(\beta'_j\circ \varphi),
\end{equation*}
where $\xi_i,\,\eta_j\in C^\infty(M_1)$ and $\alpha'_i,\,\beta'_j\in \varGamma^\infty(M_2, A_2)$, 
\begin{equation}\label{sec3:Lie algebroid morphism}
   \Phi\circ[\alpha,\,\beta]_1 = 
    \sum_{i,j}\xi_i \eta_j\, ([\alpha'_i,\,\beta'_j]_2\circ\Phi) 
    + \sum_j\bigl(\mathcal{L}_{\sharp_1(\alpha)}\eta_j\bigr) (\beta'_j\circ\Phi) 
    - \sum_i\bigl(\mathcal{L}_{\sharp_1(\beta)}\xi_i\bigr) (\alpha'_i\circ\Phi) 
\end{equation}
For further discussion of Lie algebroid morphisms, we refer to \cite{DZpoi05} and \cite{Mgen05}.

Concepts in Lie algebroid theory often appear as generalizations of standard 
notions in Poisson geometry and differential geometry. The following theorem is an 
analogue of the splitting theorem by A. Weinstein which states that any 
Poisson manifold is locally a direct product of symplectic manifold with another Poisson manifold~(see \cite{Wloc83}).
The splitting theorem for Lie algebroids appears in \cite{Dnor01, Flie02, Walm00}. 
We refer to \cite{DZpoi05} for the proof of this theorem. 

\begin{thm}\label{sec3:the splitting theorem}
{\rm (The splitting theorem \cite{Dnor01,Flie02,Walm00})}~
Let $(A\to M,\sharp,[\cdot,\,\cdot])$ be a Lie algebroid. 
For each point $m\in M$, there exist a local coordinate chart with coordinates 
$(x_1,\cdots,x_r,y_1\cdots,y_s)$ centered at $m$, where $r=\mathrm{rank}\,\sharp_p$ and $r+s=\dim M$, and 
a basis of local sections $\{\,\alpha_1,\cdots,\alpha_r,\beta_1,\cdots,\beta_s\,\}$ over an open 
neighborhood of $m$ such that 
\begin{align*}
[\alpha_j,\,\alpha_k] &= \boldsymbol{0},\quad 
[\alpha_j,\,\beta_k] = \boldsymbol{0},\quad 
[\beta_j,\,\beta_k] = \sum_\ell f^\ell_{jk}(y)\,\beta_\ell,\\
\sharp\alpha_j &= \frac{\partial}{\partial x_j},\quad 
dx_j(\sharp\beta_k) = 0\quad 
\mathcal{L}_{\frac{\partial}{\partial x_j}}\sharp\beta_k = \boldsymbol{0}
\end{align*}
for all possible indices $j,k,\ell$. Here, $f^\ell_{jk}(y)$ are smooth functions depending 
only on the variables $y=(y_1,\cdots,y_s)$. 
\end{thm}

The notion of $A$-connections given bellow generalizes 
the usual one of connections on vector bundles~(see \cite{CFint03}). 

\begin{dfn}\label{sec3:A-connection}
Let $(A\to M,\sharp,[\cdot,\,\cdot])$ be 
a Lie algebroid over $M$ and $E$ a vector bundle over $M$. 
An $\mathbb{R}$-bilinear map 
\[
 \nabla^A: 
 \varGamma^\infty(M, A)\times \varGamma^\infty(M, E)\longrightarrow \varGamma^\infty(M, E),\quad 
 (\alpha, s)\longmapsto \nabla^A_\alpha s
\]
is called an $A$-connection if it satisfies 
\begin{enumerate}[\quad \rm(1)]
 \item $\nabla^A_{f\alpha}s \;=\; f\nabla^A_\alpha s;$
 \item $\nabla^A_\alpha(fs) \;=\; f\nabla^A_\alpha s \;+\; \bigl((\sharp\alpha) f\bigr)s$ 
\end{enumerate}
for any $f\in C^\infty(M), \alpha\in\varGamma^\infty(M, A)$ and $s\in\varGamma^\infty(M, E)$. 
\end{dfn}

The notion of ordinary connection is the case where $A$ is the tangent algebroid $TM$. 
We denote by $\nabla^0$ an ordinary connection, that is, 
\[
\nabla^0: \mathfrak{X}\,(M)\times \varGamma^\infty(M, E)\longrightarrow 
\varGamma^\infty(M, E),\quad 
 (X, s)\longmapsto \nabla^0_X s. 
\]
When $E\to M$ is a complex vector bundle, an $A$-connection on $E$ 
is defined as an $A$-connection which is $\mathbb{C}$-linear on $\varGamma^\infty(M, E)$. 

Similarly to the case of usual connection theory on vector bundles, one can define 
the curvature of an $A$-connection. 
The curvature $R_\nabla^A$ of an $A$-connection $\nabla^A$ is the map 
\[
 R_\nabla^A: \varGamma^\infty(M, A)\times\varGamma^\infty(M, A)\to 
 \mathrm{End}_\mathbb{R}~(\varGamma^\infty(M, E)) 
\]
given by the usual formula 
\[
 R^A_\nabla(\alpha, \beta) = \nabla^A_\alpha\circ\nabla^A_\beta - \nabla^A_\beta\circ\nabla^A_\alpha 
 - \nabla_{[\alpha,\, \beta]}^A
\]
for any $\alpha,\,\beta\in\varGamma^\infty(M, A)$. 

For each $k\in\mathbb{N}\cup\{0\}$, 
consider the exterior bundle $\wedge^kA^*$ over $M$. 
A smooth section of $\wedge^kA^*$ is called an $A$-differential $\ell$-form. 
One defines a multilinear map, called a $A$-exterior derivative, 
$d_A:\varGamma^\infty(M, \wedge^\ell A^*)\to \varGamma^\infty(M, \wedge^{\ell +1}A^*)$ as 
\begin{align*}
 (d_A\theta)\,(\alpha_1,\cdots,\alpha_{\ell +1}) &= \sum_{j=1}^{\ell +1}(-1)^{j+1}\sharp\alpha_j~
   \bigl(\theta\,(\alpha_1,\cdots,\widehat{\alpha_j},\cdots,\alpha_{\ell +1})\bigr)\\
   &\qquad\qquad + \sum_{j<k}(-1)^{j+k}
  \theta\,\bigl([\alpha_j,\alpha_k],\alpha_1\cdots, 
    \widehat{\alpha}_j,\cdots,\widehat{\alpha}_k,\cdots,\alpha_{\ell +1}\bigr)  
\end{align*}
for any $\alpha_1,\cdots,\alpha_{k+1}\in\varGamma^\infty(M, A)$. 

\noindent The following proposition can be verified by a direct computation. 

\begin{prop}\label{sec2:prop The}
 The differential operator $d_A$ has the following properties:
 \begin{enumerate}[\quad\rm(1)]
  \item $d_A\circ d_A = 0${\rm ;}
  \item For any $A$-differential $k$-form $\theta$ and $A$-differential $\ell$-form $\vartheta$, 
      \[d_A(\theta\wedge\vartheta) \,=\, (d_A\theta)\wedge\vartheta \,+\, 
      (-1)^k\theta\wedge (d_A\vartheta). \]
 \end{enumerate}
\end{prop}
From Proposition \ref{sec2:prop The}, one finds that $(\varGamma^\infty(M, \wedge^\bullet A^*), 
d_A)$ forms a chain complex. The cohomology of $(\varGamma^\infty(M, \wedge^\bullet A^*), d_A)$ 
is called the Lie algebroid cohomology, or $A$-cohomology (see \cite{CWgeo99}). 
By the definition, the $k$-th cohomology group, denoted by 
$H^k_\mathrm{L}(M,A)$, is given by 
\[
 H^k_{\mathrm L}(M, A) = \frac{\ker\,\{d_A:\varGamma^\infty(M, \wedge^kA^*)\longrightarrow
         \varGamma^\infty(M, \wedge^{k+1}A^*)\}}
                {\mathrm{im}\,\{d_A:\varGamma^\infty(M, \wedge^{k-1}A^*)\longrightarrow
         \varGamma^\infty(M, \wedge^kA^*)\}} \, . 
\]
We denote by $[\alpha]$ the cohomology class of 
$\alpha\in\ker\,\{d_A:\varGamma^\infty(M, \wedge^kA^*)\longrightarrow
\varGamma^\infty(M, \wedge^{k+1}A^*)\}$. 
A Dirac structure $D$ over $M$ can be regarded as a Lie algebroid $D\to M$ with the bracket 
$\llbracket\cdot,\cdot\rrbracket$ and the anchor map $\sharp=\rho_{TM}={\mathrm{pr}_1}|_D$. 
The cohomology of $(M,D)$ is defined as the Lie algebroid cohomology 
$H_\mathrm{L}^\bullet(M,D)$ of the Lie algebroid $(D\to M,\,\rho_{TM},\,\llbracket\cdot,\cdot\rrbracket)$. 

Let $\phi\oplus Q$ be any $D$-differential $\ell$-form of $(M,D)$, 
where $\phi\in \Omega^\ell(M)$ and $Q\in\mathfrak{X}^\ell(M)$. 
Then, the exterior derivative $d_D(\phi\oplus Q)$ for $\phi\oplus Q$ 
is given by 
\begin{align*}
 &(d_D(\phi\oplus Q))\,\bigl(X_1\oplus \xi_1,\cdots, X_{\ell +1}\oplus \xi_{\ell +1}\bigr)\\ =&\; 
 \sum_{j=1}^{\ell +1}(-1)^{j+1}X_j~
   \Bigl(\phi\,\bigl(X_1,\cdots,\widehat{X_j},\cdots,X_{\ell +1}\bigr)
   + \bigl(\xi_1\wedge\cdots\wedge\widehat{\xi_j}\wedge\cdots\xi_{\ell +1}\bigr)\,(Q)\Bigr)\\
   &\qquad + \sum_{j<k}(-1)^{j+k}\phi\,\bigl([X_j,X_k],X_1\cdots, 
    \widehat{X}_j,\cdots,\widehat{X}_k,\cdots,X_{\ell +1}\bigr)\\ 
   &\qquad\qquad + \sum_{j<k}(-1)^{j+k}
   \Bigl(\bigl(\mathcal{L}_{X_j}\xi_k - \mathrm{i}_{X_k}d\xi_j\bigr)\wedge\xi_1\wedge\cdots\wedge 
    \widehat{\xi}_j\wedge\cdots\wedge\widehat{\xi}_k\wedge\cdots\wedge\xi_{\ell +1}\Bigr)\,(Q), 
\end{align*}
where $X_1\oplus \xi_1,\cdots, X_{\ell +1}\oplus \xi_{\ell +1}$ are smooth sections of $D$. 
Especially, if $f$ is a smooth function on $M$, $d_Df$ is calculated to be 
\begin{equation}\label{sec3:computation}
 (d_Df)(X\oplus \xi) = Xf. 
\end{equation}
As noted in the subsection 2.1, one has a skew-symmetric form 
$\Pi^\sharp:\rho_{T^*M}\,(D)\to (\rho_{T^*M}\,(D))^*$ 
by $\Pi^\sharp\,(\xi_j) = X_j$ for 
$X_j\oplus \xi_j\in\varGamma^\infty(M,D)~(\forall j=1,\cdots,\ell+1)$. 
Noting that  
\[
 \mathcal{L}_{X_j}\xi_k - \mathrm{i}_{X_k}d\xi_j = \{\,\xi_j,\, \xi_k\,\}, 
\]
we get 
\begin{equation*}
 (d_D(\phi\oplus Q))\,\bigl(X_1\oplus \xi_1,\cdots,X_{\ell +1}\oplus \xi_{\ell +1}\bigr) 
 =\; 
 (d\phi)\,(X_1,\cdots,X_{\ell +1}) + 
  (\partial Q)\,(\xi_1,\cdots,\xi_{\ell +1}), 
\end{equation*}
where $\partial:\mathfrak{X}^\bullet(M)\to \mathfrak{X}^{\bullet +1}(M)$ denotes 
the contravariant exterior derivative (see I. Vaisman \cite{Vlec94}): 
\begin{align*}
 (\partial Q)\,(\alpha_1,\cdots,\alpha_{\ell +1}) 
 &= \sum_{j=1}^{\ell +1}(-1)^{j+1}\Pi^\sharp(\alpha_j)~
   \Bigl(Q\,(\alpha_1,\cdots,\widehat{\alpha_j},\cdots,\alpha_{\ell +1})\Bigr)\\
   &\qquad\qquad + \sum_{j<k}(-1)^{j+k}Q\,\bigl(\{\alpha_j,\alpha_k\},\alpha_1\cdots, 
    \widehat{\alpha}_j,\cdots,\widehat{\alpha}_k,\cdots,\alpha_{\ell +1}\bigr),  
\end{align*}
for any $\alpha_1,\cdots,\alpha_{\ell +1}\in \Omega^1(M)$. As a result, we get the following lemma:

\begin{lem}\label{sec3:decomposition}
Let $(M,D)$ be a Dirac manifold. 
The $D$-exterior derivative 
$d_D:\varGamma^\infty(M,\,\wedge^\bullet D^*)\to \varGamma^\infty(M,\,\wedge^{\bullet +1} D^*)$ 
has the decomposition of exterior differentials $d$ and $\partial$: 
\[
 d_D(\phi\oplus Q) = d\phi + \partial Q . 
\]
\end{lem}

The anchor map $\rho_{TM}:D\to TM$ 
has the natural extension to a map $\wedge^2\rho_{TM}:\varGamma^\infty(M, \wedge^2D) 
\to \mathfrak{X}^2\,(M)$ by 
\[
 (\alpha_1\wedge \alpha_2)\,(\wedge^2\rho_{TM}\,(\vartheta)) := \rho_{TM}^*(\alpha_1)\wedge\rho_{TM}^*(\alpha_2)\,
 (\vartheta) 
\]
for any $\alpha_1, \alpha_2\in\Omega^1(M)$ and $\vartheta \in \varGamma^\infty(M,\,\wedge^2D)$. 
The dual map 
$(\wedge^2\rho_{TM})^*:\Omega^2(M)\to \varGamma^\infty(M, \wedge^2D^*)$
of $\wedge^2\rho_{TM}$ is explicitly given by 
\begin{align*}
\bigl((\wedge^2\rho_{TM})^*\sigma\bigr)\,(\psi_1,\,\psi_2) 
&= \sigma\,\bigl(\wedge^2\rho_{TM}\,(\psi_1\wedge\psi_2)\bigr)\\
&= \sigma\,\bigl(\rho_{TM}(\psi_1),\rho_{TM}(\psi_2)\bigr) 
\end{align*}
for any pair of sections $\psi_1, \psi_2\in \varGamma^\infty(M,D)$ and any $\sigma\in \Omega^2(M)$. 
Since the de Rham cohomology group $H_\mathrm{dR}^\bullet(M)$ of $M$ is isomorphic to the $\check{\rm C}$ech cohomology group 
$H^\bullet(M,\mathbb{R})$, one gets a homomorphism 
\begin{equation}\label{sec3:homomorphism}
(\wedge^2\rho_{TM})^*: H^2(M,\mathbb{R})\longrightarrow H_\mathrm{L}^2(M,D), \qquad [\sigma]
 \longmapsto [(\wedge^2\rho_{TM})^*\sigma]. 
\end{equation}
from the second cohomology group $H^2(M,\mathbb{R})$ of $M$  
to the Lie algebroid cohomology group $H_\mathrm{L}^2(M, D)$ of $(M,D)$. 
It is easy to check that the composition map of $(\wedge^2\rho_{TM})^*$ and the exterior derivative $d:\Omega^1(M)\to \Omega^2(M)$ 
commutes with the map $d_D\circ {\rho_{TM}}^*:\Omega^1(M)\to \varGamma^\infty(M,\wedge^2D^*)$.
\begin{prop}\label{sec3:commute}
 $(\wedge^2\rho_{TM})^*\circ d = d_D\circ {\rho_{TM}}^*$. 
\end{prop}
\subsection{The pull-back of Lie algebroid}

Let $(A\to M,\,\sharp,\,\llbracket\cdot,\cdot\rrbracket)$ be a Lie algebroid and 
$\Phi:M'\to M$ a smooth map. 
Assume that the differential $d\Phi$ of $\Phi$ is transversal to the anchor map 
$\sharp:A\to TM$ in the sense that 
\begin{equation}\label{sec3:transversality}
\mathrm{im}\,\sharp_{\Phi(x)} \,+\, \mathrm{im}\,(d\Phi)_x \,=\, 
         T_{\Phi(x)}M,\quad (\forall x\in M').
\end{equation}
Here, $\mathrm{im}\,\sharp_{\Phi(x)}$ stands for the image of $\sharp_{\Phi(x)}$.
This assumption leads us to the following condition:
 \begin{equation}
 \mathrm{im}\,({id}_x\times\sharp_{\Phi(x)}) \,+\, T_{(x.\Phi(x))}
  \bigl(\mathrm{graph}(\Phi)\bigr)
  \,=\, T_xM'\oplus T_{\Phi(x)}M,\quad (\forall x\in M'), 
 \end{equation}
where ${id}_x$ means the identity map on $T_xM'$.
The condition ensures that the preimage 
\begin{equation}\label{pullback}
({id}\times\sharp)^{-1}T\bigl(\mathrm{graph}(\Phi)\bigr) \,=\, 
\coprod_{x\in M'}\Bigl\{\,(V;\,\alpha)\,\bigm|\, V\in T_xM',\,\alpha\in A_{\Phi(x)},\,
(d\Phi)_x(V) = \sharp_{\Phi(x)}\alpha(\Phi(x))\,\Bigr\}
\end{equation} 
is a smooth subbundle of $(\left.TM'\times A)\right|_{\mathrm{graph}(\Phi)}$ over $\mathrm{graph}(\Phi)\approx M'$. 
The vector bundle (\ref{pullback}) is a Lie algebroid 
whose anchor map is the natural projection $\mathrm{pr}_1(V;\alpha):=V$ and whose Lie bracket is given by 
\begin{align*}
 \Bigl[ \Bigl(V \,;\, \sum_jf_j\otimes \alpha_j\Bigr),\, &
 \Bigl(V' ; \sum_kg_k\otimes \beta_k\Bigr)\Bigr] \notag \\
&= \Bigl(\,[V,\, V'] \,;\, \sum_{j,k}\,f_jg_k\llbracket\alpha_j,\, \beta_k\rrbracket + 
    \sum_k\,\mathcal{L}_Vg_k\otimes \beta_k - \sum_j\,\mathcal{L}_{V'}f_j\otimes \alpha_j \,\Bigr) 
\end{align*}
for any section written in the form 
\[
 \Bigl(V \,;\, \sum_jf_j\otimes \alpha_j\Bigr)\qquad 
\bigl(\,V\in T_xM', f_j\in C^\infty(M'),\,\alpha_j\in A_{\Phi(x)},\,
 (d\Phi)_x(V) = f_j(x)\,\sharp_{\Phi(x)}{\alpha_j}(\Phi(x))\,\bigr). 
\]
This Lie algebroid is called the pull-back of Lie algebroid and denoted by $\Phi^!A$ 
(see \cite{HMalg90}). We remark that $f^!A$ has rank 
$\mathrm{rank}(\Phi^!A) = \mathrm{rank}\, A - \dim M + \dim M'$. 

Let $\Phi:M'\to (M,D)$ be a smooth map from a smooth manifold $M'$ to a Dirac manifold $(M,D)$
which satisfies the condition (\ref{sec3:transversality}). 
Given a $D$-differential $\ell$-form $\vartheta$, 
we define a $\Phi^!D$-differential $\ell~(\ell>0)$-form $\Phi^*\vartheta$, 
called the pull-back of $\vartheta$, as 
\begin{align*}
 &(\Phi^*\vartheta)_x\bigl(\,\bigl(V_1;(d\Phi)_x(V_1),\xi_1\bigr),\cdots,
 \bigl(V_\ell;(d\Phi)_x(V_\ell),\xi_\ell\bigr)\,\bigr)\\ 
 :=\; &
 \vartheta_{\Phi(x)}\bigl(\,\bigl((d\Phi)_x(V_1),\,\xi_1\bigr),\cdots,
   \bigl((d\Phi)_x(V_\ell),\,\xi_\ell\bigr)\,\bigr)
\end{align*}
for any $V_j\in T_xM'$ and $\xi_j\in T^*_{\Phi(x)}M~(j=1,2,\cdots,\ell)$. 
If $\ell=0$, the pull-back of $0$-form $f\in C^\infty(M)$ is defined as $\Phi^*f:=f\circ\Phi$. 
By using Lemma \ref{sec3:decomposition}, 
it can be easily verified that $\Phi^*$ and $d_D$ commute with each other, that is, 
\[
 \Phi^*\circ d_D = d_D\circ\Phi^*. 
\]
\vspace{0.5cm}
 
Applying Theorem \ref{sec3:the splitting theorem} to a Dirac structure $D\to M$, we find that, 
for each $p\in M$, there exist local coordinates 
$(x_1,\cdots,x_r,y_1\cdots,y_s)$ centered at $p$ ($r=\mathrm{rank}\,(\rho_{TM})_p$ and $r+s=\dim M$) and 
a basis of local sections 
\begin{equation}\label{sec4:local sections}
 \biggl(\frac{\partial}{\partial x_1}\oplus \lambda_1\biggr),\cdots,
 \biggl(\frac{\partial}{\partial x_r}\oplus \lambda_r\biggr),\, 
 (Y_1\oplus \mu_1),\cdots,(Y_s\oplus \mu_s)
\end{equation}
over an open neighborhood $W$ of $p$ which satisfy 
\begin{equation*}\label{sec3:the splitting theorem2}
Y_k = \sum_{j=1}^s h_{jk}(y)\,\frac{\partial}{\partial y_j},\qquad 
\mathcal{L}_{\frac{\partial}{\partial x_j}}\lambda_k
 = \mathrm{i}_{\frac{\partial}{\partial x_k}}d\lambda_j,\qquad 
\mathcal{L}_{\frac{\partial}{\partial x_j}}\mu_k
 = \mathrm{i}_{Y_k}d\lambda_j. 
\end{equation*}
for all possible indices $j,k$, where we note that $\det\, \bigl(h_{jk}(y)\bigr)_{1\leq j,k\leq s}= 0$. 
Let us consider the pull-back of the Lie algebroid $D\to M$ along the projection 
$\mathrm{pr}_p:M\times \mathbb{R}\to M$: 
\[
D_1:={\mathrm{pr}_p}^!D = 
\coprod_{(p,t)\in M\times\mathbb{R}}\Bigl\{\,(X_p,f(p,t)\left({\partial}/{\partial t}\right)_t
 \,;\, X_p\oplus\xi_p )\, 
\bigm|\, X_p\oplus \xi_p\in D_p,\, f\in C^\infty(M\times\mathbb{R})\,\Bigr\}. 
\]
Noting that $\mathrm{rank}(D_1)=\mathrm{rank}\,D - \dim M + \dim\,(M\times\mathbb{R}) = \dim\,(M\times\mathbb{R})$, 
$D_1$ has the local basis of the smooth sections on $W\times\mathbb{R}$
\[
 \biggl(\frac{\partial}{\partial x_j}, \,\boldsymbol{0}\,;
   \, \frac{\partial}{\partial x_j}\oplus \lambda_j\biggr),\quad 
 (Y_k,\boldsymbol{0}\,;\, Y_k\oplus \mu_k),\quad 
 \biggl(\boldsymbol{0},\,\frac{\partial}{\partial t}\,;
   \, \boldsymbol{0}\oplus \boldsymbol{0}\biggr)\quad (1\leq j\leq r,\,1\leq k\leq s) 
\]
induced by (\ref{sec4:local sections}). 
We denote their dual basis 
by 
\begin{equation*}
 \gamma_1,\cdots, \gamma_r,\, 
  \delta_1,\cdots,\delta_s,\, dt \in \varGamma^\infty(W, {D_1}^*), 
\end{equation*}
that is, they are the local smooth sections of ${D_1}^*$ such that 
\begin{equation*}
\gamma_j\biggl(\frac{\partial}{\partial x_k}, \,\boldsymbol{0}\,;
   \, \frac{\partial}{\partial x_k}\oplus \lambda_k\biggr) =  
\delta_j(Y_k,\boldsymbol{0}\,;\, Y_k\oplus \mu_k) 
= \begin{cases} \,1\quad (j=k)\\ \,0\quad (j\ne k),\end{cases}\qquad 
dt\,\biggl(\boldsymbol{0},\,\frac{\partial}{\partial t}\,;
   \, \boldsymbol{0}\oplus \boldsymbol{0}\biggr) = 1 
\end{equation*}
and 
\begin{align*}
& \gamma_j(Y_k,\boldsymbol{0}\,;\, Y_k\oplus \mu_k) = 
\delta_j\,\biggl(\frac{\partial}{\partial x_k}, \,\boldsymbol{0}\,;
   \, \frac{\partial}{\partial x_k}\oplus \lambda_k\biggr) = 0, \\
& dt\,(Y_k,\boldsymbol{0}\,;\, Y_k\oplus \mu_k) = 
dt\,\biggl(\frac{\partial}{\partial x_k}, \,\boldsymbol{0}\,;
   \, \frac{\partial}{\partial x_k}\oplus \lambda_k\biggr) =0
\end{align*}
for all possible $j,k$. 
A smooth section $\alpha=(\,X,f(p, t)\,\partial/\partial t ; X\oplus \xi\,)\in
\varGamma^\infty(W\times\mathbb{R},D_1)$ 
is written in the form 
\[
 \left(\sum_{j,k} \biggl(u_j\frac{\partial}{\partial x_j} + v_kY_k\biggr), 
 \,f(p,t)\,\frac{\partial}{\partial t}\,;
   \, \sum_{j,k} \biggl(u_j\frac{\partial}{\partial x_j} + v_kY_k\biggr)\oplus 
   \sum_k^s \bigl(u_k\lambda_k + v_k\mu_k \bigr)\right). 
\]
Then, by a simple computation, one finds that 
\[
u_j = \gamma_j(\alpha),\quad 
v_k=\delta_k(\alpha),
\quad f(p,t)=dt\,(\alpha)
\qquad (1\leq j\leq r,\, 1\leq k\leq s). 
\]
Accordingly, from (\ref{sec3:computation}), the $D_1$-exterior derivative of a smooth function $F$ on $M\times\mathbb{R}$ 
is represented as 
\[
d_{D_1}F = \sum_{j}^r \,\frac{\partial F}{\partial x_j}\,\gamma_j 
   + \sum_k^s\,(Y_kF)\,\delta_k + \frac{\partial F}{\partial t}\,dt. 
\]

\begin{prop}\label{sec4:Poincare lemma}
 Let $(M,D)$ be a Dirac manifold and $D_1\to M$ the pull-back of 
 Lie algebroid $D$ along the natural projection 
 $\mathrm{pr}_p:M\times \mathbb{R}\to M$. Then, $\mathrm{pr}_p$ induces the isomorphism 
\[
 {\mathrm{pr}}_p^* : H_\mathrm{L}^\bullet(M, D) \xrightarrow{\cong} 
 H_\mathrm{L}^{\bullet}(M\times\mathbb{R}, D_1). 
\]
The inverse is the homomorphism 
$\iota^*:H_\mathrm{L}^{\bullet}(M\times\mathbb{R},D_1)\to H_\mathrm{L}^\bullet(M,D)$ induced from the 
inclusion map 
$\iota:M\to M\times\mathbb{R},\;p\mapsto (p,0)$. 
\end{prop}

\noindent {\em Proof.}~ 
Since $\mathrm{pr}_p\circ\iota = id$, it holds that $\iota^*\circ\mathrm{pr}_p^*=id$. 
Therefore, it is sufficient to show $\mathrm{pr}_p^*\circ\iota^*=id$ on 
$H_L^{\bullet}(M\times\mathbb{R},D_1)$ for the proof. 
For simplicity, 
we may assume that $M$ is an Euclidean space $\mathbb{R}^{\dim M}$ and 
$W$ is a star-shaped open set with respect to the origin $\boldsymbol{0}\in \mathbb{R}^{\dim M}$. 
We remark that 
any $D_1$-differential $\ell$-form $\omega$ can be written in the form 
\[
\omega = \sum_{I,I'}f_{I,I'}(p,t)\,\gamma_I\wedge\delta_{I'} 
  + \sum_{J,J'}g_{J,J'}(p,t)\,dt\wedge\gamma_J\wedge\delta_{J'}, 
\]
where $I,I'$ and $J,J'$ run over all sequences with 
$1\leq i_1<i_2<\cdots<i_c \leq r,\,1\leq i'_1<i'_2<\cdots<i'_{c'} \leq s~(c+c'=\ell)$ and 
$1\leq j_1<j_2<\cdots<j_d \leq r,\, 1\leq j'_1<j'_2<\cdots<j'_{d'} \leq s~(d+d'=\ell-1)$, 
respectively. 

For each $\ell$. we define an operator $S_\ell$ 
from $\varGamma^\infty(W\times\mathbb{R},\wedge^\ell {D_1}^*)$ to 
$\varGamma^\infty(W\times\mathbb{R},\wedge^{\ell - 1} {D_1}^*)$ 
as 
\[
 S_\ell(\omega) := \sum_{J,J'}\left(\int_0^tg_{J,J'}(p,t)\,dt\right)\gamma_J\wedge\delta_{J'}. 
\]
Then, by Proposition \ref{sec2:prop The}, we have 
\begin{align*}
 &d_{D_1}(S_\ell(\omega)) \\
 =& \sum_{J,J'}\left\{\sum_j\frac{\partial}{\partial x_j}
   \biggl(\int_0^t\,g_{J,J'}(p,t)\,dt\biggr)\,\gamma_j + 
   \sum_k Y_k\biggl(\int_0^t\,g_{J,J'}(p,t)dt\biggr)\,\delta_k + 
 \frac{\partial}{\partial t}\biggl(\int_0^t\,g_{J,J'}(p,t)\,dt\biggr)\,dt\right\}
 \gamma_J\wedge\delta_{J'}\\
 &\quad - \sum_{J,J'}\biggl(\int_0^t\,g_{J,J'}(p,t)\,dt\biggr)\,(d_{D_1}\gamma_J)
 \wedge \delta_{J'} + \sum_{J,J'}\biggl(\int_0^t\,g_{J,J'}(p,t)\,dt\biggr)\,\gamma_J
 \wedge (d_{D_1}\delta_{J'}). 
\end{align*}
On the other hand, the $D_1$-exterior derivative of $\omega$ is calculated to be 
\begin{align*}
 d_{D_1}(\omega) 
 =& \sum_{I,I'}\left\{\sum_j \frac{\partial f_{I,I'}}{\partial x_j}(p,t)\,
   \gamma_j + \sum_k (Y_k f_{I,I'})\,\delta_k + 
 \frac{\partial f_{I,I'}}{\partial t}(p,t)\,dt \right\}\gamma_I\wedge\delta_{I'}\\
 &\quad - \sum_{I,I'}f_{I,I'}(p,t)\,(d_{D_1}\gamma_I)\wedge \delta_{I'} + 
 \sum_{I,I'}f_{I,I'}(p,t)\,\gamma_I \wedge (d_{D_1}\delta_{I'})\\
 &\qquad + \sum_{J,J'}\left\{\sum_j \frac{\partial g_{J,J'}}{\partial x_j}(p,t)\,
   \gamma_j + \sum_k (Y_k g_{J,J'})\,\delta_k + 
 \frac{\partial g_{J,J'}}{\partial t}(p,t)\,dt \right\}dt\wedge\gamma_J\wedge\delta_{J'}\\
 &\qquad\quad + \sum_{J,J'}g_{J,J'}(p,t)\,dt\wedge (d_{D_1}\gamma_J)\wedge \delta_{J'}  
 - \sum_{J,J'}f_{J,J'}(p,t)\,dt\wedge\gamma_J \wedge (d_{D_1}\delta_{J'}).  
\end{align*}
Therefore,
\begin{align*}
 S_\ell(d_{D_1}(\omega)) &= 
  \sum_{I,I'}\biggl(\int_0^t\,\frac{\partial f_{I,I'}}{\partial t}(p,t)\,dt\biggr)\,
 \gamma_I\wedge\delta_{I'} 
 - \sum_{J,J'}\left\{\sum_j\biggl(\int_0^t\,\frac{\partial g_{J,J'}}{\partial x_j}(p,t)\,dt\biggr)
 \gamma_j\right\}\wedge\gamma_J\wedge\delta_{J'}\\
 &\quad - \sum_{J,J'}\left\{\sum_k\biggl(\int_0^t\,(Y_k g_{J,J'})\,dt\biggr)\,
 \delta_k\right\}\wedge\gamma_J\wedge\delta_{J'} 
 + \sum_{J,J'}\biggl(\int_0^t\,g_{J,J'}(p,t)\,dt\biggr)\,(d_{D_1}\gamma_J)\wedge\delta_{J'}\\
 &\qquad - \sum_{J,J'}\biggl(\int_0^t\,g_{J,J'}(p,t)\,dt\biggr)\,\gamma_J\wedge(d_{D_1}\delta_{J'}). 
\end{align*}
As a result, we have that 
\begin{align}\label{sec4:Poincare1}
 &d_{D_1}(S_\ell(\omega)) + S_\ell(d_{D_1}(\omega)) \notag \\
 =&\; 
  \sum_{I,I'}\biggl(\int_0^t\,\frac{\partial f_{I,I'}}{\partial t}(p,t)\,dt\biggr)\,
 \gamma_I\wedge\delta_{I'} 
 + \frac{\partial}{\partial t}\biggl(\int_0^t\,g_{J,J'}(p,t)\,dt\biggr)\,dt\wedge 
 \gamma_J\wedge\delta_{J'} \notag\\
 =&\; 
  \sum_{I,I'}f_{I,I'}(p,t)\,\gamma_I\wedge\delta_{I'} 
 - \sum_{I,I'}f_{I,I'}(p,0)\,\gamma_I\wedge\delta_{I'} 
 + \sum_{J,J'}g_{J,J'}(p,t)\, dt\wedge \gamma_J\wedge\delta_{J'} \notag \\ 
 =&\; \omega
  - \sum_{I,I'}f_{I,I'}(p,0)\,\gamma_I\wedge\delta_{I'}. 
\end{align}
Here, we recall again that 
the pull-backs 
${\mathrm{pr}}_p^*:\varGamma^\infty(M,\wedge^{\ell} D^*)\to 
\varGamma^\infty(M\times\mathbb{R},\wedge^{\ell} {D_1}^*)$ and 
$\iota^*:\varGamma^\infty(M\times\mathbb{R},\wedge^{\ell} {D_1}^*)\to 
\varGamma^\infty(M,\wedge^{\ell} D^*)$
are given by 
\begin{align*}
&({\mathrm{pr}_p}^*\vartheta)\,
\left(\biggl(X_1,f_1\,\frac{\partial}{\partial t} : X_1\oplus \xi_1\biggr),\cdots, 
 \biggl(X_\ell,f_\ell\,\frac{\partial}{\partial t} : X_\ell\oplus \xi_\ell\biggr)\right) 
 := \vartheta\,\Bigl(X_1\oplus \xi_1,\cdots,X_\ell\oplus \xi_\ell \Bigr),\\
&{\mathrm{pr}_p}^*f = f\circ{\mathrm{pr}_p}\quad (f\in C^\infty(M))
\end{align*}
and 
\begin{align*}
&(\iota^*\omega)\,(X_1\oplus \xi_1,\cdots, X_\ell\oplus \xi_\ell) := 
\omega\,\Bigl((X_1,\boldsymbol{0} : X_1\oplus \xi_1),\cdots, 
 (X_\ell,\boldsymbol{0} : X_\ell\oplus \xi_\ell)\Bigr), \\
&\iota^*F = F\circ\iota\quad (F\in C^\infty(M\times\mathbb{R})) 
\end{align*}
respectively.
By a simple computation we get
\begin{equation}\label{sec4:Poincare2}
 \omega - (\mathrm{pr}_p^*\circ\iota^*)\,\omega = 
 \omega 
 - \sum_{I,I'}f_{I,I'}(p,0)\,\gamma_I\wedge\delta_{I'}. 
\end{equation}
From (\ref{sec4:Poincare1}) and (\ref{sec4:Poincare2}) it follows that 
\[
 d_{D_1}\circ S_\ell + S_\ell\circ d_{D_1} \;=\; id - (\mathrm{pr}_p^*\circ\iota^*). 
\]
Since $d_{D_1}\circ(\iota\circ\mathrm{pr}_p)^*=(\iota\circ\mathrm{pr}_p)^*\circ d_{D_1}$, 
it turns out that $\omega = d_{D_1}\bigl(\omega-(\mathrm{pr}_p^*\circ\iota^*)\,\omega\bigr)$
for any $D_1$-differential $\ell$-form $\omega$ such that $d_D\omega = 0$. 
That is, $[\omega]=0$ in $H_{\rm L}^\ell(M\times\mathbb{R},D_1)$. This completes the proof. 
\qquad\qquad $\Box$ 
\subsection{Dirac-Chern classes of complex line bundles}

We let $q_L:L\to M$ be a complex line bundle over a Dirac manifold $(M,D)$ and 
$\{(U_j,\varepsilon_j)\}_{j}$ be a family of pairs which gives local trivializations of $L$. 
That is, $\{U_j\}_j$ is an open covering of $M$ and $\varepsilon_j$ are 
nowhere vanishing smooth sections on $U_j$ such that the map 
\[
  U_j \times \mathbb{C} \longrightarrow {q_L}^{-1}(U_j),\qquad 
  (x,\, z)\longmapsto z\,\varepsilon_j(x)
\]
for each $j$ is a diffeomorphism. 

\noindent A $D$-connection 
\[
 \nabla^D : \varGamma^\infty(M, D)\times \varGamma^\infty(M, L)\to \varGamma^\infty(M, L), 
\] 
is also considered as a map from 
$\varGamma^\infty(M,L)$ to $\varGamma^\infty(M,D^*)\otimes_{C^\infty(M)}\varGamma^\infty(M,L)$ by 
\[
\varGamma^\infty(M,L)\ni 
s\longmapsto \{\psi\mapsto \nabla^D_\psi s\}\; 
\in \mathrm{Hom}_{C^\infty(M)}\,\bigl(\varGamma^\infty(M,D),\,\varGamma^\infty(M,L)\bigr). 
\]
On each $U_j$, $\nabla^D\varepsilon_j$ is written as 
\[
 \nabla^D\varepsilon_j = 2\pi \sqrt{-1}\,\sigma_j \otimes \varepsilon_j,
\]
by using a smooth section $\sigma_j \in \varGamma^\infty(U_j, D^*)$. 
Since the transition function $g_{jk}$ on $U_j\cap U_k (\ne \emptyset)$ is given by 
$g_{jk}(x):=\varepsilon_k(x)/\varepsilon_j(x)\; (x\in U_j\cap U_k)$, 
we have $\varepsilon_k(x)=g_{jk}(x)\varepsilon_j(x)$. 
It follows from a simple computation, that 

\begin{equation}\label{sec3:computation1}
 \nabla^D \varepsilon_k = (\,2\pi\sqrt{-1}\,g_{jk}\sigma_j + d_Dg_{jk})\otimes \varepsilon_j. 
\end{equation}

On the other hand,  

\begin{equation}\label{sec3:computation2}
 \nabla^D \varepsilon_k = 2\pi\sqrt{-1}\,g_{jk}\sigma_{k}\otimes \varepsilon_j. 
\end{equation}
It immediately follows from (\ref{sec3:computation1}) and (\ref{sec3:computation2}) that 

\begin{equation}\label{sec3:1-section}
 \sigma_j - \sigma_k = \frac{\sqrt{-1}}{2\pi}\, \frac{d_Dg_{jk}}{g_{jk}}. 
\end{equation}
As a result, one gets a $D$-differential 2-form $\tau$ defined on the whole of $M$ 
by $\tau = d_D\sigma_j = d_D\sigma_k\; (U_j\cap U_k\ne \emptyset)$. 
It is easy to verify that $\tau$ satisfies 
\[
 \bigl(R_\nabla^D(\psi_1,\psi_2)\bigr)(\varepsilon_j) 
= 2\pi\sqrt{-1}\,\tau\,(\psi_1,\psi_2)\varepsilon_j 
 \quad (\forall \psi_1,\psi_2\in\varGamma^\infty(M,D))
\]
for each $j$. That is, $\tau$ is the curvature 2-section of $\nabla^D$ 
(see Remark \ref{sec3:remark} below). Obviously, $\tau$ defines a second $D$-cohomology 
class $[\tau]\in H^2_\mathrm{L}(M, D)$. 

\begin{prop}\label{sec3:the 2nd cohomology class}
The cohomology class $[\tau]$ determined by the curvature 2-section $\tau$ 
does not depend on the choice of the $D$-connection $\nabla^D$. 
\end{prop}

\noindent {\em Proof.}~
Let $\nabla'$ be another $D$-connection on $L\to M$ whose curvature is $R'$ 
and $\sigma'_j$ the corresponding local sections in $\varGamma^\infty(U_j, D^*)$. 
Denoting by $\tau'$ the curvature 2-section corresponding to $R'$, we have 
\begin{equation}\label{sec3:computation3}
 \tau' - \tau = d_D\sigma'_j - d_D\sigma_j = d_D(\sigma'_j - \sigma_j) 
\end{equation}
on each $U_j$. 
We define a $\mathbb{C}$-linear map $\widehat{\nabla}$ as 
\[
 \widehat{\nabla}:\varGamma^\infty(M, L)\longrightarrow  
 \varGamma^\infty(M, D^*)\otimes_{C^\infty(M)}\varGamma^\infty(M, L), \quad 
 s\longmapsto (\nabla' - \nabla^D)\,s. 
\]
On $U_j$, the following holds 
\[
 \widehat{\nabla}_\psi\varepsilon_j = (\nabla'_\psi - \nabla^D_\psi)\,\varepsilon_j 
 = 2\pi\sqrt{-1}\,\sigma'_j(\psi)\,\varepsilon_j  - 2\pi\sqrt{-1}\,\sigma_j (\psi)\,\varepsilon_j 
 = 2\pi\sqrt{-1}\,(\sigma'_j - \sigma_j)(\psi)\,\varepsilon_j
\]
for any $\psi$ in $\varGamma^\infty(U_j, D)$. Putting $\widehat{\sigma}_j=\sigma'_j - \sigma_j$ 
for each $j$, we find that, by $(\ref{sec3:1-section})$, 
\[
\widehat{\sigma}_j - \widehat{\sigma}_k 
= (\sigma'_j - \sigma'_k) - (\sigma_j - \sigma_k) 
= \frac{\sqrt{-1}}{2\pi}\, \frac{d_Dg_{jk}}{g_{jk}} - \frac{\sqrt{-1}}{2\pi}\, \frac{d_Dg_{jk}}{g_{jk}} = 0 
\]
on $U_j\cap U_k \ne\emptyset$. 
Accordingly, there exists a $D$-differential 1-form $\widehat{\sigma}$ over the whole of 
$M$ given by $\widehat{\sigma} = \widehat{\sigma}_j = \widehat{\sigma}_k$ 
on $U_j\cap U_k\,(\ne\emptyset)$. 
Therefore, it follows from (\ref{sec3:computation3}) that 
\[
 \tau' - \tau = d_D\widehat{\sigma}_j = d_D\widehat{\sigma}. 
\]
This shows that $[\tau] = [\tau']$ in $H_\mathrm{L}^2(M,D)$. 
\qquad\qquad\qquad\qquad\qquad\qquad\qquad\qquad\qquad\qquad\qquad\qquad$\Box$

\begin{dfn}
 Let $L\to M$ be a complex line bundle over a Dirac manifold $(M,D)$ and 
 $\nabla^D$ any $D$-connection on $L$. 
 The cohomology class $[\tau]\in H_\mathrm{L}^2(M,D)$ by the $D$-differential 2-form $\tau$ which 
 corresponds to the curvature of $\nabla^D$ is called the first Dirac-Chern class of $L\to M$. 
 We denote the first Dirac-Chern class of $L$ by $c^D_1(L)$. 
\end{dfn}

We assume that the line bundle $L\to M$ has a Hermitian metric $h$. 
A $D$-connection $\nabla^D$ is called a Hermitian $D$-connection with respect to $h$ if 
\begin{equation*}
 \rho_{TM}(\psi)\,\bigl(h(s_1,s_2)\bigr) 
 = h\,(\nabla^D_\psi s_1,\, s_2) + h\,(s_1,\, \nabla^D_\psi s_2)
\end{equation*}
for any smooth section $s_1, s_2$ of $L$ and any smooth section $\psi$ of $D$. 
The following proposition can be shown in a way similar to the case of the ordinary connections 
on Hermitian line bundles~(see \cite{Kqua70}).
\begin{prop}
 The curvature 2-section $\tau$ of $\nabla^D$ is a real $D$-differential 2-form. 
\end{prop}

\begin{rmk}\label{sec3:remark}
Let $A$ be any Lie algebroid over $M$. 
In general, an $A$-connection $\nabla^A$ on a vector bundle $\pi:E\to M$ 
is considered as an $\mathbb{R}$-linear map 
$\varGamma^\infty(M,E)$ to $\varGamma^\infty(M,A^*)\otimes_{C^\infty(M)}\varGamma^\infty(M,E)$ 
satisfying the condition (2) in Definition \ref{sec3:A-connection}. 
Let $\{V_\lambda\}_\lambda$ be an open covering which gives local trivializations of $E$ 
and $s_1,\cdots,s_r~(r=\mathrm{rank}\,E)$ be smooth sections such that $s_1(p),\cdots,s_r(p)$ 
is a basis for the fiber $\pi^{-1}(p)$ for every $p\in V_\lambda$. 
One can verify that 
there exists a matrix $\theta = (\theta_{jk})$ of local sections of $A^*$ over $V_\lambda$ 
such that 
\[
 \nabla^As_k \;=\; \sum_{j}\theta_{jk}\otimes s_j \quad 
(\,\theta_{jk}\in\varGamma^\infty(V_\lambda, A^*)\,). 
\]
The matrix $\theta$ is called a connection 1-section (see \cite{Flie02}). 
In the same manner as the ordinary connection theory, 
the curvature $R_\nabla^A$ of $\nabla^A$ is written as  
\[
\bigl(R_\nabla^A(\alpha_1,\alpha_2)\bigr)(s_k) 
= \sum_{j}\kappa_{jk}(\alpha_1,\alpha_2)s_j 
 \quad (\forall \alpha_1,\alpha_2\in\varGamma^\infty(V_\lambda, A))
\]
on each $V_\lambda$, where $\kappa_{jk}\in\varGamma^\infty(V_\lambda, \wedge^2A^*)$. 
The matrix $\kappa=(\kappa_{jk})$ 
is called the curvature 2-section of $\nabla^A~(see \ \cite{Flie02})$. 
\end{rmk}

\section{Prequantization of Dirac manifolds}

\subsection{$\Omega$-compatible Poisson structures}

Let $(M, D)$ be a Dirac manifold. 
As mentioned in Section 2, $(M, D)$ has the presymplectic structure 
$\Omega^\flat:\rho_{TM}\,(D)\to \rho_{TM}\,(D)^*$ by 
(\ref{sec2:2-form}). 
We define a singular distribution $\mathcal{V}$ 
as 
\begin{equation*}
\mathcal{V} := \ker \Omega^\flat = D \cap TM. 
\end{equation*}
Here, we remark again that $D_p\cap T_pM$ is thought of as a subspace of 
either $T_pM\oplus T_p^*M$ or $T_pM$ at each $p\in M$. 
We consider a subset $\mathcal{H}$ of $D$ whose fibers $\mathcal{H}_p$ are subspaces 
of $D_p$ satisfying 
\begin{equation}\label{sec3:disjoint union} 
 \mathcal{V}_p\oplus \mathcal{H}_p = D_p\quad (\forall p\in M), 
\end{equation} 
and fix it. We denote a singular distribution $\rho_{TM}\,(\mathcal{H})\subset TM$ by $\mathcal{H}_{TM}$. 
Note that $\dim\,(\mathcal{H}_{TM})_p={\rm rank}\,(\Omega^\flat)_p$ for any $p\in M$. 
Since $\ker\Omega^\flat \cap \mathcal{H}_{TM}=\{\boldsymbol{0}\}$, it turns out that
$\mathcal{H}_{TM}$ is isomorphic to $\widehat{\mathcal{H}}_{TM}
:=\mathrm{im}\,\Omega^\flat|_{\mathcal{H}_{TM}}=\mathrm{im}\,\Omega^\flat$ by 
the restriction map 
\[
\Omega^\flat|_{\mathcal{H}_{TM}}:\mathcal{H}_{TM}\xrightarrow{\cong} \widehat{\mathcal{H}}_{TM}.
\]
We denote its inverse map 
$(\Omega^\flat|_{\mathcal{H}_{TM}})^{-1}:\widehat{\mathcal{H}}_{TM}\to \mathcal{H}_{TM}$ by $\Theta^\sharp$. 
Then, it can be easily verified that
\begin{equation}\label{sec3:compatibility}
 \Theta^\sharp\circ \Omega^\flat|_{\mathcal{H}_{TM}} = \mathrm{id}_{\mathcal{H}_{TM}} \quad\text{\rm and}\quad 
 \Omega^\flat|_{\mathcal{H}_{TM}}\circ \Theta^\sharp = \mathrm{id}_{\widehat{\mathcal{H}}_{TM}}. 
\end{equation}
As mentioned in Section 2, there exists a bundle map $\Pi_p^\sharp$ defined as 
\[
 {\mathcal{V}_p}^\circ\ni\eta_p\longmapsto 
 \bigl\{\,\xi_p\mapsto \,\xi_p\,\bigl(\Pi_p^\sharp(\eta_p)\bigr):=\xi_p\,(Y_p)\,\bigr\}
 \in T_pM/(D_p\cap T_pM). 
\]
We here remark that ${\mathcal{V}_p}^\circ$ is the annihilator of $\mathcal{V}_p$ in $T_p^*M$. 
From the definition of $\Omega$ and Proposition \ref{sec2:prop}, 
the image $\mathrm{im}\,\Omega^\flat$ of $\Omega^\flat$ 
turns out to be 
\[
\mathrm{im}\,\Omega_p^\flat = \rho_{T^*M}\,(D_p) = (D_p\cap T_pM)^\circ = {\mathcal{V}_p}^\circ. 
\]
So, we have $\widehat{\mathcal{H}}_{TM}= \mathcal{V}^\circ$, and 
find that $\Theta^\sharp = \Pi^\sharp|_{\widehat{\mathcal{H}}_{TM}}$. 
By Proposition \ref{sec2:prop}, any admissible function $f\in C^\infty_{\rm adm}(M,D)$ satisfies 
\[
df \in \widehat{\mathcal{H}}_{TM}. 
\] 
This allows us to define a vector field 
\[
H_f:=\Theta^\sharp(df)\in \mathcal{H}_{TM}. 
\]
Since $f$ is admissible, there exists a vector field $X_f$ such that 
$X_f\oplus df\in \varGamma^\infty(M, D)$. 
It is easy to see that $((H_f)_p-(X_f)_p)\oplus \boldsymbol{0}\in \mathcal{V}_p\subset D_p$
at each $m\in M$. It follows from this that 
\[
(H_f)_p\oplus (df)_p = \bigl((H_f)_p-(X_f)_p\bigr)\oplus \boldsymbol{0} 
+ (X_f)_p\oplus (df)_p \in D_p. 
\]
So, it turns out that $H_f\oplus df\in \varGamma^\infty(M,D)$. 
For any $f,g\in C^\infty_{\rm adm}(M,D)$, let us define their bracket $\{f,g\}$ as 

\begin{equation}\label{sec3:bracket}
\{f,\,g\}:= H_gf. 
\end{equation}
It is easily verified that, for any $f,g\in C_\mathrm{adm}^\infty(M,D)$, 
\[
 \{f,g\}=\Omega\,(H_f, H_g). 
\]
Since $H_f\oplus df,\, H_g\oplus dg\in\varGamma^\infty(M, D)$, we have that 
\[
 [H_g,H_f]\oplus d\{f,g\} = \llbracket\,H_g\oplus dg,\,H_f\oplus df\,\rrbracket\in \varGamma^\infty(M, D).
\]
So, $d\{f,g\}$, also, is the admissible function. 
This implies that one can define the operator 
\[
\{\cdot,\cdot\}:C^\infty_{\rm adm}(M,D)\times C^\infty_{\rm adm}(M,D)
\longrightarrow C^\infty_{\rm adm}(M,D), 
\]
as (\ref{sec3:bracket}), which is both bilinear and skew-symmetric. 
Furthermore, it turns out that this bracket $\{\cdot,\,\cdot\}$ coincides with $\{\cdot,\,\cdot\}'$ defined by (\ref{sec2:bracket}) 
since $(H_g-X_g)f=0$ for $(H_g-X_g)\oplus \boldsymbol{0}\in \mathcal{V}$. Consequently, we obtain the following 
proposition:

\begin{prop}\label{sec3:bracket2}
 $(\,C^\infty_{\rm adm}(M,D),\,\{\cdot,\cdot\}\,)$ forms a Poisson algebra, which does not depend on the choice of 
 $\mathcal{H}$. Moreover, it holds that 
 \begin{equation}\label{sec4:formula}
  [H_f, H_g] + H_{\{f,g\}} = \boldsymbol{0}
 \end{equation}
 for any $f,g\in C^\infty_{\rm adm}(M,D)$. 
\end{prop}
Following \cite{Vgeo83}, we say that 
the bracket $\{\cdot,\,\cdot\}$ by (\ref{sec3:bracket}) is an $\Omega$-compatible Poisson 
structure. 

\begin{ex}
Let us consider a Dirac manifold $(\mathbb{R}^4, \mathrm{graph}\,(\omega^\flat))$ by 
the presymplectic form in Example \ref{sec2:example presymplectic}. 
The presymplectic form $\Omega$ is entirely $\omega$ and written in the matrix form
\[
 \Omega^\flat = \omega^\flat = 
 \left(\begin{array}{rrrr} 0&-1&0&-1\\1&0&0&0\\0&0&0&0\\1&0&0&0 \end{array}\right). 
\]
From this, one can write the Dirac structure as 
\[
 \mathrm{graph}\,(\omega^\sharp) = \mathrm{span}\,\left\{\,
 \frac{\partial}{\partial x_1}\oplus (dx_2+dx_4),\, \frac{\partial}{\partial x_2} \oplus (-dx_1),\, 
 \frac{\partial}{\partial x_3}\oplus \boldsymbol{0},\, \frac{\partial}{\partial x_4}\oplus (-dx_1)\right\}. 
\]
The subspace $\mathcal{V}= \mathrm{graph}\,(\omega^\flat)\cap T\mathbb{R}^4$ 
is given by 
\[
  \mathcal{V} = \mathrm{span}\,\left\{\, \frac{\partial}{\partial x_2}-\frac{\partial}{\partial x_4},\, 
 \frac{\partial}{\partial x_3}\,\right\}. 
\]
Then, we can take a subspace $\mathcal{H}$ as 
\[
 \mathcal{H} = \mathrm{span}\,\left\{\, 
 \frac{\partial}{\partial x_1}\oplus (dx_2+dx_4),\, \frac{\partial}{\partial x_2} \oplus (-dx_1) \right\}. 
\]
One easily checks that $\mathrm{graph}\,(\omega_p^\flat) = \mathcal{V}_p\oplus \mathcal{H}_p$. 
The vector field $H_f$ for the admissible function 
$f(\boldsymbol{x})={x_1}^2 + k\,(x_2 + x_4)~(k\in\mathbb{R})$ is given by 
\[
H_f =  k\,\frac{\partial}{\partial x_1} - 2x_1\,\frac{\partial}{\partial x_2}. 
\]
\end{ex}

\begin{ex}
Let $(M, F\oplus F^\circ)$ be a Dirac manifold obtained from an involutive distribution $F\subset TM$ of constant rank~
(see Example \ref{sec2:regular foliation}). Since any vector field $X\in F$ is embedded 
in $F\oplus F^\circ$ with $X\hookrightarrow (X,\boldsymbol{0})$, one finds that 
$\mathcal{V}=F\oplus\{\boldsymbol{0}\}\cong F$. It follows from this that 
$\mathcal{H}=\{\boldsymbol{0}\}\oplus F^\circ\cong F^\circ$. 
If $f\in C^\infty_{\rm adm}(M,F\oplus F^\circ)$, the vector field $H_f$ for f is given by 
$H_f=\boldsymbol{0}$. 
\end{ex}

\begin{ex}
Consider the Dirac manifold $(\mathbb{R}^2,\mathrm{graph}\,(\pi^\sharp))$ 
induced by a Poisson bivector $\pi=G(\boldsymbol{x}) \partial/\partial x_1
\wedge \partial/\partial x_2$, where $G(\boldsymbol{x})=G(x_1,x_2)$ is a smooth function on 
$\mathbb{R}^2$. We remark that the Dirac structure $\mathrm{graph}\,(\pi^\sharp)$ is 
written in the form 
\[
 \mathrm{graph}\,(\pi^\sharp) = \mathrm{span}\,\left\{\, 
 G(\boldsymbol{x})\,\frac{\partial}{\partial x_1}\oplus dx_2,\, 
 - G(\boldsymbol{x})\,\frac{\partial}{\partial x_2}\oplus dx_1 \,\right\} 
\]
and any smooth function on $(\mathbb{R}^2,\mathrm{graph}\,(\pi^\sharp))$ is admissible. 
The distribution $\mathcal{V}$ is given by $\mathcal{V}=\mathrm{graph}\,(\pi^\sharp)\cap 
T\mathbb{R}^2= \{\boldsymbol{0}\}$ and consequently, 
$\mathcal{H}$ is $\mathcal{H}=\mathrm{graph}\,(\pi^\sharp)$. 
For a smooth function $h$, the vector field $H_h$ is represented as 
\[
 H_h = G(\boldsymbol{x})\,\left(\frac{\partial h}{\partial x_2}\frac{\partial}{\partial x_1} 
 - \frac{\partial h}{\partial x_1}\frac{\partial}{\partial x_2}\right). 
\]
\end{ex}

\begin{ex}\label{sec4:example,contact manifolds} 
Let $D_{\eta}$ be a Dirac structure obtained by a contact manifold $(\mathbb{R}^{2n+1},\, \eta=dz-\sum_{i=1}^n y_idx_i)$ 
(see Example \ref{sec2:example,1-form}). It is easily checked that $\mathcal{V}$ is the subbundle generated by 
the Reeb vector field $\frac{\partial}{\partial z}$. Accordingly, $\mathcal{H}$ can be taken as 
\[
\mathcal{H}=\left\{\,\sum_{i=1}^n a_i\,\Bigl(\frac{\partial}{\partial x_i}\oplus dy_i\Bigr) + 
 \sum_{j=1}^n b_j\,\Bigl(\frac{\partial}{\partial y_j}\oplus (-dx_j)\Bigr)\,\Bigm|\, a_i,\,b_j\in C^\infty(\mathbb{R}^{2n+1})~
 (\forall i,j=1,\cdots,n)\,\right\}. 
\] 
Then, the vector field $H_f$ for $f\in C^\infty(\mathbb{R}^{2n+1})$ is represented as 
\[
 H_f = \sum_{i=1}^n\left(\frac{\partial f}{\partial y_i}\frac{\partial}{\partial x_i}-
               \frac{\partial f}{\partial x_i}\frac{\partial}{\partial y_i}\right). 
\]
\end{ex}

\begin{ex}
Consider a Dirac structure $\mathrm{graph}\,(\omega)\subset\mathbb{T}M$ over $M=\mathbb{R}^3$ 
induced from a presymplectic form $\omega=F\, dx_1\wedge dx_2 + G\, dx_2\wedge dx_3 + H\, dx_3\wedge dx_1$, where 
$F,\,G$ and $H$ are smooth functions on $M$ such that 
$\frac{\partial F}{\partial x_3}(p) + \frac{\partial G}{\partial x_1}(p) + \frac{\partial H}{\partial x_2}(p)=0~(\forall p\in M)$. That is, 
\[
 \mathrm{graph}\,(\omega) = \mathrm{span}\,\Bigl\{\,\frac{\partial}{\partial x_1}\oplus (F\,dx_2-H\,dx_3),\,
      \frac{\partial}{\partial x_2}\oplus (G\,dx_3-F\,dx_1),\,\frac{\partial}{\partial x_3}\oplus (H\,dx_1-G\,dx_2)\,\Bigr\}. 
\] 
At each point $p=(x_1,x_2,x_3)$ where $F(p)\ne 0,\,G(p)\ne 0,\,H(p)\ne 0$, the distribution $\mathcal{V}$ is given by 
$\mathcal{V} = \mathrm{span}\,\{G\frac{\partial}{\partial x_1}+H\frac{\partial}{\partial x_2}+F\frac{\partial}{\partial x_3}\}$. 
Accordingly, we can take $\mathcal{H}$ as the subspace generated by 
\[
 \Bigl(\,F\,\frac{\partial}{\partial x_2}-H\,\frac{\partial}{\partial x_3}\,\Bigr) \oplus 
  \left\{\frac{GH}{F}\,(Fdx_2-Hdx_3)+\frac{F^2+H^2}{F}\,(Gdx_3-Fdx_1)\right\} 
\]
and
\[
  \Bigl(\,G\,\frac{\partial}{\partial x_3} - F\,\frac{\partial}{\partial x_1}\,\Bigr)\oplus 
  \left\{-\frac{F^2+G^2}{F}\,(Fdx_2-Hdx_3)-\frac{GH}{F}\,(Gdx_3-Fdx_1)\right\}.
\]
Then, the inverse map $\Theta^\sharp$ of $\Omega^\flat$ restricted to $\mathcal{H}_{TM}$ is represented in the matrix form 
\[
\Theta^\sharp = \frac{1}{F(F^2+G^2+H^2)}\left(\begin{array}{@{\,}cc@{\,}}
  -GH&F^2+G^2\\ -(F^2+H^2)& GH\end{array}\right). 
\]
with respect to the basis $\{\,F\,dx_2-H\,dx_3,\,G\,dx_3-F\,dx_1\,\}\subset \widehat{\mathcal{H}}_{TM}$
and $\Bigl\{\,F\,\frac{\partial}{\partial x_2}-H\,\frac{\partial}{\partial x_3},\,
  G\,\frac{\partial}{\partial x_3}-F\,\frac{\partial}{\partial x_1}\,\Bigr\}\subset \mathcal{H}_{TM}$. 
If $f$ is an admissible function, $f$ shall satisfy the condition that $G\,\frac{\partial f}{\partial x_1}+H\,\frac{\partial f}{\partial x_2}
  +F\,\frac{\partial f}{\partial x_3}=0$. 
Then, the vector field $H_f$ for $f$ is given by 
\begin{align*}
 H_f &= -\frac{1}{F^2\,(F^2+G^2+H^2)}\,\left\{\,(F^2+G^2)\,\frac{\partial f}{\partial x_1}+GH\,\frac{\partial f}{\partial x_2}\,\right\}\,
 \Bigl(\,F\,\frac{\partial}{\partial x_2}-H\,\frac{\partial}{\partial x_3}\,\Bigr)\\ &\qquad\qquad  
  -\frac{1}{F^2\,(F^2+G^2+H^2)}\,\left\{\,GH\,\frac{\partial f}{\partial x_1}+(F^2+H^2)\,\frac{\partial f}{\partial x_2}\,\right\}\,
\Bigl(\,G\,\frac{\partial}{\partial x_3}-F\,\frac{\partial}{\partial x_1}\,\Bigr). 
\end{align*}
On the other hand, at each point $p=(x_1,\,x_2,\,x_3)$ where $F(p)\ne 0,\,G(p)\ne 0,\,H(p)= 0$, $\mathcal{V}$ is spanned 
by $G\,\frac{\partial}{\partial x_1}+F\,\frac{\partial}{\partial x_3}$. Define $\mathcal{H}$ as 
\[
\mathcal{H} = \mathrm{span}\,\left\{\,\frac{\partial}{\partial x_2}\oplus (G\,dx_3-F\,dx_1),\quad 
\Bigl(G\,\frac{\partial}{\partial x_3}-F\,\frac{\partial}{\partial x_1}\Bigr)\oplus \{-(F^2+G^2)\,dx_2\}\,\right\}, 
\]
and we find that $\Theta^\sharp$ is given by 
\[
\Theta^\sharp(dx_2) = -\frac{1}{F^2+G^2}\,\Bigl(G\,\frac{\partial}{\partial x_3}-F\,\frac{\partial}{\partial x_1}\,\Bigr),\qquad
\Theta^\sharp(G\,dx_3-F\,dx_1) = \frac{\partial}{\partial x_2}. 
\]
Therefore, the vector field $H_f$ for a function $f$ with $G\,\frac{\partial f}{\partial x_1}+F\,\frac{\partial f}{\partial x_3}=0$ 
is represented as 
\[
 H_f = -\frac{1}{F}\,\frac{\partial f}{\partial x_1}\,\frac{\partial}{\partial x_2} 
 -\frac{1}{F^2+G^2}\frac{\partial f}{\partial x_2}\,\Bigl(G\,\frac{\partial}{\partial x_3}-F\,\frac{\partial}{\partial x_1}\,\Bigr) 
\]
Lastly, at each point $p=(x_1,\,x_2,\,x_3)$ where $F(p)\ne 0,\,G(p)= 0,\,H(p)= 0$, 
$\mathcal{V}=\mathrm{span}\,\Bigl\{\frac{\partial}{\partial x_3}\Bigr\}$. 
Define $\mathcal{H}$ as 
\[
\mathcal{H} = \mathrm{span}\,\left\{\,\frac{\partial}{\partial x_1}\oplus F\,dx_2,\quad 
\frac{\partial}{\partial x_2}\oplus (-F\,dx_1)\,\right\}, 
\]
and we find that 
\[
\Theta^\sharp (dx_1) = -\frac{1}{F}\,\frac{\partial}{\partial x_2},\qquad 
\Theta^\sharp (dx_2) = \frac{1}{F}\,\frac{\partial}{\partial x_1}. 
\]
If $f$ is a function which satisfies $\frac{\partial f}{\partial x_3}=0$, the vector field $H_f$ is given by 
\[
 H_f = \frac{1}{F}\,\left(\,\frac{\partial f}{\partial x_2}\,\frac{\partial}{\partial x_1} 
        - \frac{\partial f}{\partial x_1}\,\frac{\partial}{\partial x_2}\,\right).
\]
\end{ex}

\subsection{Quantizable Dirac manifolds}

We let $(M, D)$ be a Dirac manifold and fix a singular distribution $\mathcal{H}\subset D$ for $\mathcal{V}=\ker\,\Omega$. 
Suppose that there exists a line bundle $L\overset{q_L}{\to} M$ over $(M, D)$ 
with a $D$-connection $\nabla^D$ whose curvature is $R_\nabla^D$. 
For $\mathcal{H}$, 
we define a map~ $\hat{} : C^\infty_{\rm adm}(M,D)\to \mathrm{End}_\mathbb{C}\,(\varGamma^\infty(M, L))$ 
from the Poisson algebra $(\,C^\infty_{\rm adm}(M,D),\,\{\cdot,\cdot\}\,)$ 
to Lie algebra $(\mathrm{End}_\mathbb{C}\,(\varGamma^\infty(M, L)),\,[\cdot, \cdot])$ as 
\begin{equation}\label{sec4:prequantization}
 \hat{f}s := -\nabla_{H_f\oplus df}^Ds - 2\pi\sqrt{-1}fs \qquad (\forall s\in\varGamma^\infty(M, L))
\end{equation}
for each $f\in C^\infty_{\rm adm}(M, D)$. 

\begin{prop}\label{sec4:prequantization2}
 The map~ $\hat{} : C_{\rm adm}^\infty(M, D)\to \mathrm{End}_\mathbb{C}(\varGamma^\infty(M, L))$ 
 is a representation of $C^\infty_{\rm adm}(M, D)$ on $\varGamma^\infty(M, L)$, that is, it holds that 
 \begin{equation}\label{sec4:bracket preserving}
  \widehat{\{f,g\}} = [\hat{f},\,\hat{g}]
 \end{equation}
 for all $f,g\in C^\infty_{\rm adm}(M, D)$ if and only if 
 \begin{equation}\label{sec3:prequantization}
  R^D_\nabla\bigl(H_f\oplus df, H_g\oplus dg\bigr) = 2\pi\sqrt{-1}\,\Lambda( H_f\oplus df,\,H_g\oplus dg ), 
 \end{equation}
where $\Lambda$ is the skew-symmetric pairing 
$\Lambda(\cdot,\cdot) := \langle\cdot, \cdot\rangle_-$ in Section 2. 
\end{prop}

\noindent{\em Proof}.~ 
Using Proposition \ref{sec3:bracket2}, we have that 
\begin{align*}
 [\hat{f},\,\hat{g}]s &= 
  \hat{f}(\hat{g}s) - \hat{g}(\hat{f}s)\\
 &= \hat{f}\,(-\nabla^D_{H_g\oplus~ dg}s - 2\pi\sqrt{-1}gs) - 
  \hat{g}\,(-\nabla^D_{H_f\oplus~ df}s - 2\pi\sqrt{-1}fs)\\
 &= \nabla^D_{H_f\oplus~ df}\circ\nabla^D_{H_g\oplus~ dg}s - \nabla^D_{H_g\oplus~ dg}\circ\nabla^D_{H_f\oplus~ df}s\\ 
 &\qquad\qquad\qquad 
 - 2\pi\sqrt{-1}\,\Bigl\{g\nabla^D_{H_f\oplus~ df}s - f\nabla^D_{H_g\oplus~ dg}s + \nabla^D_{H_g\oplus~ dg}(fs) 
       - \nabla^D_{H_f\oplus~ df}(gs)\Bigr\}\\
 &= R_\nabla^D\bigl(H_f\oplus df,\, H_g\oplus dg\bigr)s - \nabla^D_{H_{\{f,g\}}\oplus~ d\{f,g\}}s - 
 4\pi\sqrt{-1}\,\{f,g\}s\\
 &= \widehat{\{f,g\}}s + R_\nabla^D\bigl(H_f\oplus df,\,H_g\oplus dg\bigr)s - 2\pi\sqrt{-1}\,\{f,g\}s 
\end{align*}
for any admissible function $f, g$ on $(M, D)$ and any smooth section $s$ of $L\to M$.  
The bracket $\{f,g\}$ is calculated to be 
\begin{align*}
 \{f, g\} &= \frac{1}{2}\,\bigl(\{f,g\}-\{g,f\}\bigr)
   = \frac{1}{2}\,\bigl(df(H_g) - dg(H_f)\bigr)\\
  &= \langle\,H_f\oplus df\, H_g\oplus dg\,\rangle_-. 
\end{align*}
From this, we immediately get (\ref{sec3:prequantization}) as 
the necessary and sufficient condition for
the map~ $\hat{}$~ to preserve their brackets. 
\qquad\qquad\qquad\qquad\qquad\qquad\qquad\qquad\qquad\qquad\qquad\qquad\qquad
\qquad\qquad\qquad\qquad\quad $\Box$

\begin{dfn}
A Dirac manifold $(M,D)$ is said to be prequantizable if there exists a Hermitian 
line bundle $(L,h)$ over $M$ 
with a Hermitian $D$-connection $\nabla^D$ in the sense that 
\begin{equation}\label{sec4:hermitian}
H_f\,\bigl(h(s_1,s_2)\bigr) 
 = h\,\bigl({\nabla}^D_{H_f\oplus~ df} s_1,\, s_2\bigr) 
 + h\,\bigl(s_1,\, \nabla^D_{H_f\oplus~ df} s_2\bigr), 
\end{equation}
which satisfies the condition 
{\rm (\ref{sec3:prequantization})}. The line bundle is called the prequantum bundle. 
\end{dfn}

Let us consider the skew-symmetric pairing $\Lambda$ again. We find that 
$\Lambda:\varGamma^\infty(M, D)\times\varGamma^\infty(M, D)\to C^\infty(M)$ is closed 
with regard to the differential operator $d_D$. 
Indeed, by Lemma \ref{sec2:integrability} and the Cartan formula, $d_D\Lambda$ is calculated to be 
\begin{align*}
 &(d_D\Lambda)\bigl(X\oplus \xi,\,Y\oplus \eta,\,Z\oplus \zeta\bigr)\\
 =&\; X\,\bigl(\eta(Z)-\zeta(Y)\bigr) - Y\,\bigl(\xi(Z)-\zeta(X)\bigr) 
    + Z\,\bigl(\xi(Z)-\eta(X)\bigr)\\ 
  &\qquad\qquad - \Lambda\,\bigl(\llbracket X\oplus \xi,\, Y\oplus \eta\rrbracket,\, Z\oplus \zeta\bigr) 
  - \Lambda\,\bigl(\llbracket Y\oplus \eta,\,Z\oplus \zeta\rrbracket,\, X\oplus \xi \bigr) 
  - \Lambda\,\bigl(\llbracket Z\oplus \zeta,\, X\oplus \xi\rrbracket,\, Y\oplus \eta\bigr)\\
 =&\; X\,\bigl(\eta(Z)-\zeta(Y)\bigr) - Y\,\bigl(\xi(Z)-\zeta(X)\bigr) + Z\,\bigl(\xi(Z)-\eta(X)\bigr) 
  - (\mathcal{L}_X\eta)(Z) - (\mathcal{L}_Y\zeta)(X) - (\mathcal{L}_Z\xi)(Y)\\
 &\qquad\qquad  + (d\xi)(Y,Z) + (d\eta)(Z,X) + (d\zeta)(X,Y) + \xi([Y,Z]) 
   + \eta([Z,X]) + \zeta([X,Y])\\
 =&\; 0 
\end{align*}
for any section $X\oplus \xi,\, Y\oplus \eta$ and $Z\oplus \zeta$ of $D$. 
Accordingly, the $D$-differential 2-form $\Lambda$ defines the second cohomology class $[\Lambda]$ 
in the Lie algebroid cohomology. 
Additional to this, we have that 
\begin{align*}
 \Lambda(X\oplus \xi, Y\oplus \eta) &=\; \frac{1}{2}\,\bigl\{\xi\,(Y) - \eta\,(X)\bigr\}\\ 
  &=\; \frac{1}{2}\,\bigl\{\Omega\,(X, Y) - \Omega\,(Y,X)\bigr\} = \Omega\,(X, Y)\\ 
  &=\; ((\wedge^2\rho_{TM})^*\Omega)\,\bigl(X\oplus \xi,\,Y\oplus \eta\bigr). 
\end{align*}
That is, it holds that 
\begin{equation}\label{sec4:relation}
 \Lambda \;=\; (\wedge^2\rho_{TM})^*\Omega. 
\end{equation}

\begin{thm}\label{sec4:main theorem}
A Dirac manifold $(M,D)$ is prequantizable if and only if 
the $D$-cohomology class $[\Lambda]$ of $\Lambda$ lies in the image $\iota_*(H^2(M;\mathbb{Z}))${\rm :}  
\begin{equation}\label{sec4:new integrability}
[\Lambda] \in (\wedge^2\rho_{TM})^*\Bigl(\iota_*(H^2(M,\mathbb{Z}))\Bigr)\subset H_\mathrm{L}^2(M, D), 
\end{equation}
where $\iota_*$ is the map from $H^\bullet(M,\mathbb{Z})$ to $H^\bullet(M,\mathbb{R})$ induced from 
the inclusion $\iota:\mathbb{Z}\hookrightarrow\mathbb{R}$.
\end{thm}

\noindent {\em Proof}.~ 
We assume that 
\begin{equation}
[\Lambda] \in (\wedge^2\rho_{TM})^*\Bigl(\iota_*(H^2(M,\mathbb{Z}))\Bigr). 
\end{equation}
Let $\{W_j\}_j$ be a contractible open covering of $M$ and $\beta\in\Omega^2(M)$ a closed 2-form on $M$. 
By Poincar$\acute{\rm e}$'s lemma, there exist 1-forms $\alpha_j\in \Omega^1(W_j)$ such that 
\[
 \beta|_{W_j} = d\alpha_j 
\]
on each $W_j$. We here remark that $W_j\cap W_k$ is also contractible whenever $W_j$ and $W_k$ are so. 
Accordingly, by using Poincar$\acute{\rm e}$'s lemma again, one can write 
\[
 \alpha_j - \alpha_k = dw_{jk}
\]
for some function $w_{jk}\in C^\infty(W_j\cap W_k)$ on $W_j\cap W_k \ne\emptyset$. 
As a result, we obtain 
$D$-differential 1-forms $\sigma_j\in\varGamma^\infty(W_j,\,D^*)$ which satisfy 
\[
 \Lambda|_{W_j} = d_D\sigma_j 
\]
on each $W_j$ and find that the functions $\{w_{jk}\}_{j,k}$ satisfy 
\begin{equation}\label{sec4:proof}
 \sigma_j - \sigma_k = d_Dw_{jk}.
\end{equation} 
on $W_j\cap W_k\ne\emptyset$ by using Proposition \ref{sec3:commute}. 


We remark that $f_{jk\ell}:=w_{jk} + w_{k\ell} - w_{j\ell}$ are constant functions which take the values 
in $\mathbb{Z}$~(see \cite{Kqua70}) and 
define functions $c_{jk}\in C^\infty(W_j\cap W_k)$ as $c_{jk}:=\exp(-2\pi \sqrt{-1}\,w_{jk})$. 
By (\ref{sec4:proof}) , we have 
\[
  \sigma_j - \sigma_k = \frac{\sqrt{-1}}{2\pi}\,\frac{d_D\,c_{jk}}{c_{jk}}. 
\]
In addition, we can verify that the functions $\{c_{jk}\}_{j,k}$ satisfy the cocycle condition:
\[
 c_{jk}c_{k\ell} = \exp\,(-2\pi\sqrt{-1}\,(w_{jk} + w_{k\ell})) = 
 \exp\,(-2\pi\sqrt{-1}\,f_{jk\ell})\exp\,(-2\pi\sqrt{-1}\,w_{j\ell}) = c_{j\ell}
\] 
on $W_j\cap W_k\cap W_\ell\,(\ne\emptyset)$. 
Consequently, one can obtain a line bundle $L\to M$ whose transition functions are 
$\{c_{jk}\}_{j,k}$ and on which $\{\sigma_{j}\}_j$ determine a connection $\nabla^D$
 with curvature $\Lambda$. 

We define a Hermitian metric $h$ on $L$ as 
\[
 h_p\,(s_1,s_2) := \overline{z}_1z_2, 
\]
for any section $s_1(p)=(p,z_1),\,s_2(p)=(p,z_2)\in W_j\times\mathbb{C}$ on each 
open set $W_j$ of the trivialization, 
and fix $\mathcal{H}$ for the singular distribution $\mathcal{V}=\ker\Omega^\flat$. 
Then, $\nabla^D$ turns out to be a Hermitian connection in the sense of (\ref{sec4:hermitian}). 
Indeed, letting $s_1,s_2$ be smooth sections locally written in the form 
$s_1(p)=g_1(p)\,\varepsilon_j(p),\,s_2(p)=g_2(p)\,\varepsilon_j(p)~(g_1(p), g_2(p)\in\mathbb{C})$, 
where $\varepsilon$ is the nowhere 
vanishing section, 
and $f\in C^\infty_{\rm adm}(M,D)$, we have 
\begin{align*}
 \bigl(H_f\bigl(h\,(s_1,s_2)\bigr)\bigr)(p) &= 
 \bigl(H_f\bigl(h\,(g_1\,\varepsilon_j,g_2\,\varepsilon_j)\bigr)\bigr)(p)
 = \bigl(H_f(\overline{g}_1g_2)\bigr)(p)\\ 
 &=(H_f\overline{g}_1)(p)\,g_2(p) + \overline{g_1(p)}(H_fg_2)(p). 
\end{align*}
Recall again that $H_f$ is the vector field which is uniquely determined for $f$ with respect to $\mathcal{H}$. 

On the other hand, 
\begin{align*}
 &h\,\bigl({\nabla}^D_{H_f\oplus~df} s_1,\, s_2\bigr)(p) 
 + h\,\bigl(s_1,\, \nabla^D_{H_f\oplus~df} s_2\bigr)(p)\\
=&\, h\,\Bigl((H_fg_1)\varepsilon_j 
   + 2\pi\sqrt{-1}\, g_1\sigma_j(H_f)\varepsilon_j,\, g_2\,\varepsilon_j\Bigr)(p) 
 + h\,\Bigl(g_1\varepsilon_j,\,(H_fg_2)\varepsilon_j + 2\pi\sqrt{-1}\, g_2\sigma_j(H_f)\varepsilon_j\Bigr)(p)\\
=&\, \overline{(H_fg_1)_p}g_2(p) - 2\pi\sqrt{-1}~\overline{g_1(p)\sigma_j(H_f\oplus df)}g_2(p) 
 + \overline{g_1(p)}(H_fg_2)_p + 2\pi\sqrt{-1}~\overline{g_1(p)}g_2(p)\sigma_j(H_f\oplus df). 
\end{align*}
From the assumption, each connection 1-section $\sigma_j$ is real. Accordingly, 
\[
 h\,\bigl({\nabla}^D_{H_f\oplus~df} s_1,\, s_2\bigr)(p) 
 + h\,\bigl(s_1,\, \nabla^D_{H_f\oplus~df} s_2\bigr)(p) = 
 (H_f\overline{g}_1)_pg_2(p) + \overline{g_1(p)}(H_fg_2)_p. 
\]
Therefore, we have that 
\[
H_f\,\bigl(h(s_1,s_2)\bigr) 
 = h\,\bigl({\nabla}^D_{H_f\oplus~df} s_1,\, s_2\bigr) 
 + h\,\bigl(s_1,\, \nabla^D_{H_f\oplus df} s_2\bigr). 
\]
This results in that $(M,D)$ is prequantizable. 

Conversely, we suppose that $(M, D)$ is prequantizable, 
that is, there is the prequantization bundle $(L,\nabla^D)$ over $M$. 
Note that the $D$-differential 2-form which corresponds to $R^D_\nabla$ is $\Lambda$: 
\[
 R_\nabla^D = 2\pi\sqrt{-1}\,\Lambda. 
\]
It is well-known that 
the isomorphism classes of Hermitian line bundles over $M$ are classified by 
the second cohomology classes 
through the map which assigns to the isomorphism class of a line bundle $K\to M$ 
the first Chern class $c_1(K)\in H^2(M, \mathbb{Z})$ 
of $K$~(see \cite{Kqua70} and \cite{Wgeo91}).  
According to this, one obtains an ordinary Hermitian connection $\nabla^0$ on $L$ whose curvature is 
$R_0$. The curvature form  
$F_{\nabla^0}$ corresponding to $R_0$ satisfies 
\[
 c_1(L)=[F_{\nabla^0}]\in H^2(M, \mathbb{Z}).
\] 
The map $R_1:\varGamma^\infty(M,D)\times\varGamma^\infty(M,D)\to 
\mathrm{End}_\mathbb{C}\,(\varGamma^\infty(M,L))$ defined as 
\[
 R_1 := R_0\circ(\rho_{TM}\times\rho_{TM}) = (\wedge^2\rho_{TM})^*R_0 
\]
is the curvature of a $D$-connection 
$\nabla_1 := \nabla^0\circ(\rho_{TM}\times {\it id})$ on $L$. 
Then, the $D$-differential 2-form $\tau_1$ corresponding to 
$R_1$ is represented as $\tau_1 = (\wedge^2\rho_{TM})^*F_{\nabla^0}$
by using $F_{\nabla^0}$. 
Using Proposition \ref{sec3:the 2nd cohomology class}, 
we find that 
\[
 [\Lambda] = c_1^D(L) = [\tau_1] 
= [(\wedge^2\rho_{TM})^*F_{\nabla^0}]\in (\wedge^2\rho_{TM})^*\Bigl(\iota_*(H^2(M, \mathbb{Z}))\Bigr). 
\]
This completes the proof of Theorem \ref{sec4:main theorem}. 
\qquad\qquad\qquad\qquad\qquad\qquad\qquad\qquad\qquad\qquad\qquad\qquad $\Box$

\begin{rmk}
We would like to remark that the formula (\ref{sec4:prequantization}) and Proposition \ref{sec4:prequantization2} 
are quite different from Lemma 6.1 in \cite{WZvar05}, though they might look similar to it. 
The prequantization formula discussed in \cite{WZvar05} acts on local sections $s$ of 
a line bundle $L$ satisfying the condition $\nabla^D_{V\oplus \boldsymbol{0}}s=0~(\forall V\oplus\boldsymbol{0}\in\mathcal{V})$, 
that is, it is a representation of $C^\infty_{\rm adm}(M,D)$ on the space 
\[
\{\,s\,|\,\text{$s$ is a local section of $L$ such that $\nabla^D_{V\oplus \boldsymbol{0}}s=0
 ~(\forall V\oplus\boldsymbol{0}\in\mathcal{V})$}\,\}. 
\]
On the other hand, (\ref{sec4:prequantization}) acts on global sections of $L$ which do not necessarily require the condition 
and is a representation of $C^\infty_{\rm adm}(M,D)$ on $\varGamma^\infty(M,L)$. 
\end{rmk}

The quantization procedure in (\ref{sec4:prequantization}) generally depends on the choice of $\mathcal{H}$. 
In other words, there are as many quantization procedures of $(M, D)$ as there are the choice of $\mathcal{H}$. 
Suppose that we take $\mathcal{H}'$ different from $\mathcal{H}$. Let ${\hat{f}}'$ be the prequantization procedure 
for $f\in C^\infty_{\rm adm}(M,D)$ with respect to $\mathcal{H}'$. Then, we have that 
\[
(\hat{f}-{\hat{f}}')s := \nabla_{(H'_f-H_f)\oplus \boldsymbol{0}}^Ds \qquad (s\in \varGamma^\infty(M,L)). 
\]
From this formula, it turns out that ${\hat{f}}'$ coincides with $\hat{f}$ 
on the space 
\[
\{\,s\in \varGamma^\infty(M,L)\,|\,\nabla^D_{V\oplus \boldsymbol{0}}s=0~(\forall V\oplus\boldsymbol{0}\in\mathcal{V})\,\}. 
\]
If $\mathcal{V}=\ker\Omega^\flat\subset TM$ is a subbundle, then it is integrable and 
its leaf space $M_{\rm red}$ can be endowed with a Poisson structure $\Pi_{\rm red}$~(Corollary 2.6.3 in \cite{Cdir90}). 
Denote the natural projection from $M$ to $M_{\rm red}$ 
by $q_\mathcal{V}$ and define a subbundle $D_{\rm red}\subset \mathbb{T}M$ as 
\[
 (D_{\rm red})_p = \bigl\{\,(dq_\mathcal{V})_p(X)\oplus \overline{\xi} \,\bigm|\,X\in T_pM,\,
   \overline{\xi}\in T_{[p]}^*M_{\rm red},\, X\oplus (dq_\mathcal{V})_{[p]}^*\overline{\xi}\in D_p\,\bigr\},
\]
where $[p]:=q_\mathcal{V}(p)\in M_{\rm red}$. It can be shown that $D_{\rm red}$ is a Dirac structure over $M_{\rm red}$ and 
coincides with a Dirac structure $\mathrm{graph}\,(\Pi_{\rm red})$ induced from $\Pi_{\rm red}$. 
As we shall see in Example \ref{sec4:Poisson manifolds}, a Poisson manifold whose second cohomology class $[\pi]$ is 
integral has the prequantization bundle. In a word, we can obtain the following:

\begin{prop}
 Let $(M,\,D)$ be a Dirac manifold. If $\mathcal{V}$ is a subbundle of $TM$, there is a prequantization procedure 
 on the leaf space $M_{\rm red}$, which does not depend on the choice of $\mathcal{H}$. 
\end{prop}

We end the subsection with some examples. 

\begin{ex}\label{sec4:symplectic}  
(Symplectic manifolds)~
We let $(M,\omega)$ be a symplectic manifold and consider the Dirac structure 
$\mathrm{graph}\,(\omega^\flat)\subset TM\oplus T^*M$ 
induced from the symplectic form $\omega$~(see also Example \ref{sec2:example}). 
It is verified that the skew-symmetric 2-cocycle $\Omega$ by (\ref{sec2:2-form}) is 
entirely $\omega$.  
Therefore, it follows from the non-degeneracy of $\omega$ that 
\[
 \mathcal{V}_p=\ker \Omega_p^\flat = \ker\omega_p = \{\boldsymbol{0}\} \quad (\forall p\in M). 
\] 
Accordingly, we can take $\mathrm{graph}\,(\omega^\flat)$ 
as a subbundle $\mathcal{H}$ satisfying (\ref{sec3:disjoint union}). 
Then, for any smooth function $f$ on $M$, there exists a unique vector field $H_f$ 
such that $df=\omega^\flat(H_f)$. 
Therefore, the $\Omega$-compatible Poisson structure coincides with the natural 
Poisson structure induced from $\omega$. 
In this case, the skew-symmetric pairing $\Lambda$ is written as 
\[
\Lambda\,\bigl(X\oplus \omega^\flat(X),\,Y\oplus \omega^\flat(Y)\bigr) = \omega\,(X,Y) 
\]
and the integrability condition 
(\ref{sec4:new integrability}) is given by $[\omega]\in \iota_*(H^2(M,\mathbb{Z}))$. 
\end{ex}

\begin{ex}(Presymplectic manifolds)~
As discussed in Example \ref{sec2:example}, given a presymplectic manifold $(M,\omega)$, 
one obtains a Dirac manifold $(M,\,\mathrm{graph}\,(\omega^\flat))$. 
Similarly to Example \ref{sec4:symplectic}, $\Omega$ is entirely $\omega$. 
The singular distribution $\mathcal{V}=\ker\omega^\flat$ is given by 
$\mathcal{V}=\{\,X\in TM\,|\,\omega^\flat(X)=0\,\}\subset \mathrm{graph}\,(\omega^\flat)$. 
We let $\mathcal{H}\subset \mathrm{graph}\,(\omega^\flat)$ be a distribution which satisfies (\ref{sec3:disjoint union}) 
and fix it. 
If $f$ is a function such that there exists a vector field $X$ which satisfies $\mathrm{i}_X\omega = df$, 
one can get a unique vector field $H_f$ which belongs to $\mathcal{H}$ and  
a Poisson algebra $(\,C^\infty_\mathrm{adm}(M, \mathrm{graph}\,(\omega^\flat)), \{\cdot,\, \cdot\}\,)$. 
According to Theorem \ref{sec4:main theorem}, $(M,\,\mathrm{graph}\,(\omega^\flat))$ is prequantizable as a Dirac manifold if and only if 
the 2-form $\omega$ satisfies $[\omega]\in \iota_*(H^2(M, \mathbb{Z}))$. 
\end{ex}

\begin{ex}\label{sec4:Poisson manifolds}(Poisson manifolds)~
Let us consider the case of a Poisson manifold $(P,\pi)$. 
As seen in Example \ref{sec2:example2}, $(P, \pi)$ defines a Dirac manifold 
$(P,\,\mathrm{graph}\,(\pi^\sharp))$. 
One easily finds that 
$\mathcal{V}$ is given by $\mathcal{V}=\mathrm{graph}\,(\pi^\sharp)\cap TP = \{\boldsymbol{0}\}$. 
This permits us to take a subbundle 
$\mathrm{graph}\,(\pi^\sharp) = \{\,\pi^\sharp(\alpha)\oplus \alpha\,|\,\alpha\in T^*P\,\}$ as $\mathcal{H}$. 
Obviously, every smooth function is admissible function. 
The skew-symmetric pairing $\Lambda$ is written as 
\[
\Lambda\,\bigl(\pi^\sharp(\alpha)\oplus \alpha,\,\pi^\sharp(\beta)\oplus \beta\bigr) 
 = \pi\,(\alpha,\beta). 
\]
Then, the integrability condition 
(\ref{sec4:new integrability}) indicates that $[\pi]\in H_\mathrm{LP}^2(P,\,\pi)$ is an integral 
cohomology class, where $H_\mathrm{LP}^\bullet(P,\,\pi)$ denotes the Lichnerowicz-Poisson cohomology~(see \cite{DZpoi05, Vlec94}). 
Recall that $[\pi]$ is said to be integral if there exists a closed 2-form $\beta$ such that $[\beta]\in \iota_*(H^2(P, \mathbb{Z}))$ 
which satisfies $[\wedge^2\pi^\sharp(\beta)] = [\pi]$. 
Here, $\wedge^2\pi^\sharp$ denotes a map from $\Omega^2(P)$ to $\mathfrak{X}^2(P)$ by 
\[
 \bigl((\wedge^2\pi^\sharp)\phi\bigr)(\alpha_1,\,\alpha_2) := \phi\,(\pi^\sharp(\alpha_1),\pi^\sharp(\alpha_2)) \qquad 
 (\alpha_1,\,\alpha_2\in \Omega^1(P))
\]
Consequently, $(P,\,\pi)$ is prequantizable 
as a Poisson manifold if and only if $(P,\,\mathrm{graph}\,(\pi^\sharp))$ is prequantizable as a Dirac manifold. 
\end{ex}

\begin{ex}(Contact manifolds)~Let $M$ be a contact manifold with a contact 1-form $\eta$. As seen in Example \ref{sec2:example,1-form}, 
a subbundle $D_\eta=\mathrm{graph}\,(d\eta^\flat)$ defines a Dirac structure over $M$. In this case, the integrability condition 
(\ref{sec4:new integrability}) is represented as $[d\eta]\in \iota_*(H^2(M, \mathbb{Z}))$. 
Especially, a standard contact manifold $(\mathbb{R}^{2n+1},\,dz-\sum_{i=1}^n y_idx_i)$ turns out to be prequantizable 
as a Dirac manifold by $d\eta = \sum_i dx_i\wedge dy_i$. 
\end{ex}

\begin{ex}
Let us consider a Dirac manifold $(T^*M\times\mathbb{R},\,D_{\omega_0})$~(see Example \ref{sec2:example, canonical}). 
For any section $\bigl(X,\,f\,\frac{d}{dt}\bigr) \oplus \mathrm{pr}_1^*(\mathrm{i}_{X}\omega_0),\,
\bigl(Y,\,g\,\frac{d}{dt}\bigr) \oplus \mathrm{pr}_1^*(\mathrm{i}_{Y}\omega_0)$ of $D_{\omega_0}$, 
$\Lambda$ is written as 
\[
 \Lambda\Bigl(\bigl(X,\,f\,\frac{d}{dt}\bigr) \oplus \mathrm{pr}_1^*(\mathrm{i}_{X}\omega_0),\,
 \bigl(Y,\,g\,\frac{d}{dt}\bigr) \oplus \mathrm{pr}_1^*(\mathrm{i}_{Y}\omega_0)\Bigr) = \omega_0(X,\,Y). 
\]
Note again that $\omega_0$ is the canonical symplectic form on $T^*M$. 
Hence, the integrability condition (\ref{sec4:new integrability}) is equivalent to $[\omega_0]\in 
\iota_*(H^2(T^*M, \mathbb{Z}))$. 
From this, it follows that the Dirac manifold $(T^*M\times\mathbb{R},\,D_{\omega_0})$ is prequantizable. 
\end{ex}

\begin{ex}(Almost cosymplectic manifolds)~ Let $(M,\,D_{\omega,\eta})$ be a Dirac manifold over an almost cosymplectic manifold 
$(M,\,\omega,\,\eta)$ with $d\omega=0$ discussed in Example \ref{sec2:example, almost cosymplectic}. 
A 2-form $\Omega$ associated to $D_{\omega,\eta}$ is $\omega+d\eta$. Hence, the singular distribution $\mathcal{V}$ is given by 
$\mathcal{V}=\{\,X\in TM\,|\,\mathrm{i}_X(\omega+d\eta)=0\,\}$. Take a singular distribution $\mathcal{H}$ so that 
$\mathcal{V}$ and $\mathcal{H}$ satisfy (\ref{sec3:disjoint union}). Then, similarly to the case of presymplectic manifolds, 
one can obtain a unique vector field belonging to $\mathcal{H}$ and a Poisson algebra on the space of smooth functions $f$ 
such that $\mathrm{i}_X(\omega+d\eta)=df$ for some vector field $X$ on $M$. 
The Dirac manifold $(M,\,D_{\omega,\eta})$ is prequantizable if and only if $\omega$ and $\eta$ satisfy 
$[\omega + d\eta]\in \iota_*(H^2(M, \mathbb{Z}))$. 
\end{ex}

\section{Quantization}
\subsection{$\alpha$-density bundles} 

Let $V$ be an $n$-dimensional vector space over $\mathbb{C}$ and $\alpha$ a positive number. 
A function 
$\kappa:V\times\cdots\times V\rightarrow \mathbb{C}$
on $n$-copies of $V$ is called an $\alpha$-density over $V$ if it satisfies 
\[
\kappa\,(Av_1,\cdots,Av_n) = |\det\,A|^\alpha\,\kappa(v_1,\cdots,v_n) \quad (v_1,\cdots,v_n\in V)
\]
for any invertible linear transformation $A:V\to V$. 
Denoting the set of all $\alpha$-densities over $V$ by $\mathrm{H}^{(\alpha)}(V)$, we can check easily that 
$\mathrm{H}^{(\alpha)}(V)$ is a vector space over $\mathbb{C}$. 
Since $A\in\mathrm{GL}(V)$ acts transitively on bases in $V$, an $\alpha$-density is 
determined by its value on a single basis. 
For an alternating covariant $n$-tensor $\omega$, the map 
$|\omega|^{(\alpha)}:V\times\cdots\times V\to \mathbb{C}$ defined as 
\[
 |\omega|^{(\alpha)}\,(v_1,\cdots,v_n) 
  := |\omega\,(v_1,\cdots,v_n)|^{\alpha}\quad (v_1,\cdots,v_n\in V) 
\]
is an $\alpha$-density over $V$. If $\omega$ is nonzero, $\mathrm{H}^{(\alpha)}(V)$ is a 1-dimensional 
vector space spanned by $|\omega|^{(\alpha)}$. So, any element $\kappa\in \mathrm{H}^{\alpha}(V)$ is 
represented as $\kappa = c|\omega|^{(\alpha)}$ for some $c\in\mathbb{C}$. 

Let $M$ be a smooth manifold and $\alpha$ a positive number. The vector bundle 
\[
 {\rm H}^{\alpha} := \coprod_{m\in M}\,\mathrm{H}^{(\alpha)}(T_pM)
\]
over $M$ is called the $\alpha$-density bundle of $M$. Especially, $\mathrm{H}^{1/2}$ is called the 
half-density bundle. 
Let $(U_\lambda;(x_\lambda^1,\cdots,x_\lambda^n))$  
be local coordinate chart on $M$ and $\omega_\lambda=dx_\lambda^1\wedge\cdots\wedge dx_\lambda^n$. 
Then, a local trivialization on $U_\lambda$ is defined to be the map 
\[
\Phi_\lambda:\pi^{-1}(U_\lambda) \longrightarrow U_\lambda \times\mathbb{C},\quad 
\Phi_\lambda\Bigl(z\,|(\omega_\lambda)_p|\Bigr) = (p,z). 
\]
Letting $(U_\mu;(x_\mu^1,\cdots,x_\mu^n))$ be another chart with 
$U_\lambda\cap U_\mu\ne\emptyset$ and $\omega_\mu=dx_\mu^1\wedge\cdots\wedge dx_\mu^n$, 
we have that 
\begin{align*}
  \Phi_\lambda\circ\Phi_\mu^{-1}(p,z) &= \Phi_\lambda\Bigl(z\,|(\omega_\mu)_p|\Bigr) 
 = \Phi_\lambda
    \Biggl(
     z\left|\det\left(\frac{\partial x_\mu^k}{\partial x_\lambda^j}\right)\right|^{\alpha} 
    \!|(\omega_\lambda)_p|\Biggr)\\
&= \Biggl(p,\,z
  \left|\det\left(\frac{\partial x_\mu^k}{\partial x_\lambda^j}\right)\right|^{\alpha}\,\Biggr).     
\end{align*}
That is, ${\rm H}^{\alpha}$ is a complex 
line bundle whose transition functions are the square roots of the absolute values of 
the determinants of the matrices 
\[
\left(\frac{\partial x_\mu^k}{\partial x_\lambda^j}\,(p)\right)_{1\leq j,\,k\leq n}
\qquad (n=\dim M)
\]
by the coordinate transformations 
$x_\lambda^j = x_\lambda^j(x_\mu^1,\cdots,x_\mu^n)~(j=1,\cdots, n)$. 
A section of ${\rm H}^{\alpha}$ is called an $\alpha$-density over $M$. 
When $\alpha=1/2$, a section of ${\rm H}^{1/2}$ is called the half-density. 
As in the linear case, any $\alpha$-density $\kappa$ over $U$ can be written in the form 
$\kappa = f|\omega|^{(\alpha)}$ 
for some complex-valued function $f$. It is easily verified that 
$\mathrm{H}^\alpha\otimes\mathrm{H}^\beta\cong\mathrm{H}^{\alpha+\beta}$. Accordingly, for any half densities 
$\kappa_1,\kappa_2$ over $M$, we get a 1-density $\kappa_1\otimes\kappa_2$. 

Suppose that $(U,\phi)$ is a local coordinate chart on $M$ and $\kappa$ is a half density over $M$ 
such that the support $\mathrm{supp}\,\kappa$ of $\kappa$ is contained in $U$. 
The integral of $\kappa$ over $M$ is defined as 
\[
 \int_p\,\kappa := \int_{\phi^{-1}(U)}(\phi^{-1})^*\kappa. 
\]
We here remark that the right-hand side is represented as 
\[
 \int_{\phi^{-1}(U)}(\phi^{-1})^*\kappa = \int_{\phi^{-1}(U)}\,
  f\,|dx_1\wedge\cdots\wedge dx_n|^{(1/2)}. 
\]
If $\kappa$ is any density over $M$, the integral of $\kappa$ over $M$ is defined as 
\[
 \int_p\,\kappa := \sum_j\int_p\,\varrho_j\kappa, 
\]
where $\{\varrho_j\}_j$ means a partition of unity subordinate to smooth atlas of $M$. 

\subsection{Polarization}

The pairings $\langle\cdot,\,\cdot\rangle_\pm$ and the bracket $\llbracket\cdot,\cdot\rrbracket$ on 
$\varGamma^\infty(M,\mathbb{T}M)$ are 
naturally extended to operations on the space $\varGamma^\infty(M, \mathbb{T}M\otimes\mathbb{C})$ by 
\begin{align*}
&\bigl\langle(X\oplus\xi) + \sqrt{-1}\,(X'\oplus\xi'),\, 
   (Y\oplus \eta)+\sqrt{-1}\, (Y'\oplus \eta')\bigr\rangle_\pm \\
&= \frac{1}{2}\,\Bigl\{\,
\xi(Y) \pm \eta(X) - \xi'(Y') \mp \eta'(X') + \sqrt{-1}\,\bigl(
\xi(Y') \pm \eta(X') + \xi'(Y) \pm \eta'(X)\bigr)\,
\Bigr\}
\end{align*}
and 
\begin{align*}
&\llbracket\,(X\oplus\xi) + \sqrt{-1}\,(X'\oplus\xi'),\,
   (Y\oplus \eta)+\sqrt{-1}\, (Y'\oplus \eta')\,\rrbracket\\
&= \Bigl(\,[X,Y]-[X',Y']+\sqrt{-1}\,\bigl([X',Y] + [X,Y']\bigr)\Bigr) \\
&\qquad\qquad
\oplus \Bigl(\mathcal{L}_X\eta - \mathrm{i}_Yd\xi - \mathcal{L}_{X'}\eta' + \mathrm{i}_{Y'}d\xi' 
+ \sqrt{-1}\,
\bigl(\mathcal{L}_{X}\eta' - \mathrm{i}_Yd\xi' + \mathcal{L}_{X'}\eta - \mathrm{i}_{Y'}d\xi\bigr)\,
\Bigr), 
\end{align*}
respectively. If a complex subbundle $\mathcal{D}\subset \mathbb{T}M\otimes\mathbb{C}$ whose sections are closed under 
the extended bracket $\llbracket\cdot,\,\cdot\rrbracket$ is maximally isotropic with respect to 
the extended symmetric pairing $\langle\cdot,\,\cdot\rangle_+$, $\mathcal{D}$ is called a complex Dirac structure. 
Let $(M,D)$ be a Dirac manifold and $\Omega$ the 2-cocycle associated to it. 
Then, the complexification $D^\mathbb{C}:=D\otimes_\mathbb{R}\mathbb{C}$ of $D$ turns out to be the complex Dirac structure. 
We introduce the notion of polarization for Dirac manifold as follows: 
\begin{dfn}\label{sec5:dfn polarization}
A complex subbundle 
$\mathcal{P}\subset D^\mathbb{C}$ is called 
a (complex) polarization of $(M,\,D)$ if it satisfies the following conditions:
\begin{enumerate}[\quad\rm(P1)]
 \item $\mathcal{V}_p\otimes_\mathbb{R}\mathbb{C}\subset \mathcal{P}_p\quad (\forall p\in M)\mathrm{;}$
 \item $\Lambda\,(\tilde{X}\oplus\tilde{\xi},\,\tilde{Y}\oplus \tilde{\eta})= 0
 \quad\bigl(\forall \tilde{X}\oplus \tilde{\xi},\, \tilde{Y}\oplus \tilde{\eta} 
       \in\varGamma^\infty(M, \mathcal{P})\bigr)\mathrm{;}$
 \item $\llbracket\,\varGamma^\infty(M, \mathcal{P}), 
       \varGamma^\infty(M, \mathcal{P})\,\rrbracket 
        \subset \varGamma^\infty(M, \mathcal{P})$.
\end{enumerate}
\end{dfn}
The condition (P1) can be written in the explicit form 
\[
 \Omega(X,Y) - \Omega(X',Y') + \sqrt{-1}\,\bigl(\Omega(X',Y) + \Omega(X,Y')\bigr) = 0, 
\]
where $\tilde{X}\oplus \tilde{\xi} = (X+\sqrt{-1}X')\oplus (\Omega^\flat(X) + \sqrt{-1}\,\Omega^\flat(X'))$ and 
$\tilde{Y}\oplus \tilde{\eta} =(Y+\sqrt{-1}Y')\oplus (\Omega^\flat(Y) + \sqrt{-1}\,\Omega^\flat(Y'))$. 

\begin{ex}
 Let us consider a Dirac manifold $(\mathbb{R}^{2n},\,\mathrm{graph}\,(\omega^\flat))$ induced by 
 the standard symplectic form $\omega = \sum_j dq_j\wedge dp_j$. 
 Then, a complex subbundle $\mathcal{P}_1$ by 
 \[
   \mathcal{P}_1=\mathrm{span}\,\biggl\{\,\frac{\partial}{\partial q_j}\oplus dp_j\,
   \bigm|\, j=1,\cdots, n\,\biggr\}
 \]
 defines a polarization of $(\mathbb{R}^{2n},\,\mathrm{graph}\,(\omega^\flat))$. On the other hand, 
  \[
   \mathcal{P}_2=\mathrm{span}\,\biggl\{\,\frac{\partial}{\partial p_j}\oplus dq_j\,
   \bigm|\, j=1,\cdots, n\,\biggr\},
 \]
 also, is a polarization of $(\mathbb{R}^{2n},\,\mathrm{graph}\,(\omega^\flat))$.
\end{ex}

\begin{ex}
Consider a Poisson structure $\pi=x_3\,\frac{\partial}{\partial x_1}\wedge\frac{\partial}{\partial x_2}+
x_1\,\frac{\partial}{\partial x_2}\wedge\frac{\partial}{\partial x_3} + x_2\,\frac{\partial}{\partial x_3}
\wedge\frac{\partial}{\partial x_1}$ 
on $\mathbb{R}^3$. A complex subbundle $\mathcal{P}$ of $\mathrm{graph}\,(\pi^\sharp)^\mathbb{C}$ by 
\[
 \mathcal{P}=\mathrm{span}\,\biggl\{\,\biggl\{\Bigl(x_2\,\frac{\partial}{\partial x_3}-x_3\,\frac{\partial}{\partial x_2}\Bigr) + 
 \sqrt{-1}\,\Bigl(x_3\,\frac{\partial}{\partial x_1}-x_1\,\frac{\partial}{\partial x_3}\Bigr)\biggr\}\oplus (dx_1 + \sqrt{-1}\,dx_2)\,\biggr\}. 
\]
 turns out to be a polarization of $(\mathbb{R}^3,\,\mathrm{graph}\,(\pi^\sharp))$ by a simple computation. 
\end{ex}

\begin{ex}
 Let $(M,\omega,J)$ be a K$\ddot{a}$hler manifold, where $J$ denotes an almost complex structure on $M$.
 Define a complex subbundle $\mathcal{P}$ of $\mathrm{graph}\,(\omega^\flat)^\mathbb{C}$ as 
 \[
  \mathcal{P} = \{\,(X+\sqrt{-1}JX)\oplus (\mathrm{i}_X\omega + \sqrt{-1}\,\mathrm{i}_{JX}\omega)\,|\,X\in \mathfrak{X}(M)\,\}. 
 \]
 By a simple computation, we can check that the conditions {\rm (P1)} and {\rm (P2)} in Definition \ref{sec5:dfn polarization} are satisfied. 
 From the fact that $J$ is integrable, if follows that {\rm (P3)} holds. Therefore, the subbundle $\mathcal{P}$ is a polarization 
 of $(M,\,\mathrm{graph}\,(\omega^\flat))$. 
\end{ex}

\begin{ex}
 Let us consider a Dirac manifold $(\mathbb{R}^3,\,D_\eta)$ from a contact 1-form $\eta=dx_3 - x_2\,dx_1$. 
 As discussed in Example \ref{sec4:example,contact manifolds}, $\mathcal{V}=\ker\Omega^\flat$ is given by 
 $\mathcal{V}=\mathrm{span}\,\{\frac{\partial}{\partial x_3}\}$. A complex subbundle $\mathcal{P}$ of ${D_\eta}^\mathbb{C}$ by 
 \[
 \mathcal{P}=\mathrm{span}\,\biggl\{\,\Bigl(\,\frac{\partial}{\partial x_1}-\frac{\partial}{\partial x_2}\Bigr) \oplus  
 (dx_1+dx_2),\, \frac{\partial}{\partial x_3}\oplus \boldsymbol{0}\,\biggr\} 
\] 
 satisfies the conditions {\rm (P1)} -- {\rm (P3)} in Definition \ref{sec5:dfn polarization}. Hence, $\mathcal{P}$ is a polarization 
 of $(\mathbb{R}^3,\,D_\eta)$. 
\end{ex}

As discussed in the previous section, for a singular distribution $D\cap TM\subset D$, 
we choose a subset $\mathcal{H}\subset D$ which satisfies (\ref{sec3:disjoint union}) and fix it. 
Suppose that $(M,D)$ is prequantizable and endowed with a polarization $\mathcal{P}$. 
Let $L\to M$ be its prequantum bundle with 
the Hermitian metric $h$ on $L$ and $\nabla^D$ a Hermitian $D$-connection 
with respect to $h$.  
We remark that $\nabla^D:\varGamma^\infty(M, D)\to \mathrm{End}_\mathbb{C}(\varGamma^\infty(M, L))$ 
has a natural extension to a map 
$\nabla^D:\varGamma^\infty(M, D^\mathbb{C})\to \mathrm{End}_\mathbb{C}(\varGamma^\infty(M, L))$. 
Using the extension $\nabla^D$, 
we define a map 
\[
 \delta:\varGamma^\infty(M, D^\mathbb{C})\times \varGamma^\infty(M, L\otimes {\rm H}^{1/2})
 \to \varGamma^\infty(M, L\otimes {\rm H}^{1/2})
\] 
as 
\[
 \delta_\psi\,(s\otimes \kappa) := (\nabla^D_\psi\,s)\otimes \kappa 
 + s\otimes (\mathcal{L}_{\rho_{TM}(\psi)}\kappa), 
\]
for any $\psi\in\varGamma^\infty(M, D^\mathbb{C}),\,s\otimes\kappa
  \in\varGamma^\infty(M, L\otimes {\rm H}^{1/2})$. 
It is easily verified that $\delta$ is a $D$-connection on $L\otimes {\rm H}^{1/2}$. 
Then,  
the representation (\ref{sec4:prequantization}) of $C^\infty_{\rm adm}(M,D)$ 
can be extended to a map 
\[
 \hat{} : C^\infty_{\rm adm}(M,D)\to \mathrm{End}_\mathbb{C}\,
 (\varGamma^\infty(M, L\otimes {\rm H}^{1/2}))
\] 
by setting 
\[
 \hat{f}(s\otimes \kappa) = -\delta_{H_f\oplus~df}(s\otimes \kappa) 
    - 2\pi\,\sqrt{-1}f\,(s\otimes\kappa). 
\]
The $\hat{f}$ is also represented as 
\begin{equation}\label{sec5:representation}
 \hat{f}(s\otimes \kappa) = (\hat{f}s)\otimes\kappa - s\otimes(\mathcal{L}_{H_f}\kappa).
\end{equation}
Since $(M,D)$ is prequantizable, 
we can check that 
\[
 \widehat{\{f,g\}}\,(s\otimes\kappa) = [\hat{f},\,\hat{g}]\,(s\otimes\kappa)
\]
for all $f,g\in C^\infty_{\rm adm}\,(M,D)$ 
in the same manner as the proof for Proposition \ref{sec4:prequantization2}. 
\begin{lem}\label{sec5:lemma}
It holds that 
\[
 \delta_\psi\,\bigl(\hat{f}\,(s\otimes\kappa)\bigr)
  = \hat{f}\,\bigl(\delta_\psi(s\otimes\kappa)\bigr) 
 - \delta_{\llbracket\psi,\,H_f\oplus df\rrbracket}(s\otimes\kappa) 
\]
for any $\psi\in\varGamma^\infty(M, D^\mathbb{C}),\,s\otimes\kappa
\in\varGamma^\infty(M, L\otimes {\rm H}^{1/2})$ 
and $f\in C^\infty_{\rm adm}\,(M,D)$. 
\end{lem}

\noindent {\em Proof.}~ For any smooth section $\psi=(X,\xi)$ of $D^\mathbb{C}$, $s\otimes \kappa$ of 
$L\otimes {\rm H}^{1/2}$ and any admissible function $f$, we have that 
\begin{align*}
 \delta_\psi\circ\hat{f}\,(s\otimes\kappa) &= 
  - \delta_\psi\,\Bigl(\delta_{H_f\oplus~ df}(s\otimes \kappa) + 2\pi\,\sqrt{-1}f\,(s\otimes\kappa)\Bigr)\\ 
 &= -\delta_\psi\,\Bigl((\nabla^D_{H_f\oplus~ df}\,s)\otimes \kappa 
  + s\otimes \mathcal{L}_{H_f}\kappa\Bigr) 
  - 2\pi\,\sqrt{-1}\,\delta_\psi\bigl(f(s\otimes\kappa)\bigr)\\
 &= - (\nabla^D_\psi\circ\nabla^D_{H_f\oplus~ df}\,s)\,\otimes\kappa 
  - (\nabla^D_{H_f\oplus~ df}s)\otimes \mathcal{L}_X\kappa 
  - (\nabla^D_\psi s)\otimes \mathcal{L}_{H_f}\kappa \\
 &\qquad - s\otimes (\mathcal{L}_X\circ\mathcal{L}_{H_f}\kappa) 
  - 2\pi\,\sqrt{-1}\,(Xf)s\otimes\kappa - 2\pi\,\sqrt{-1} f\,((\nabla^D_\psi\, s)\otimes\kappa) 
  - 2\pi\,\sqrt{-1}f\,(s\otimes \mathcal{L}_X\kappa). 
\end{align*}
On the other hand, 
\begin{align*}
 \hat{f}\circ\delta_\psi(s\otimes\kappa) &= 
  \hat{f}\,\Bigl((\nabla^D_\psi s)\otimes\kappa + s\otimes\mathcal{L}_X\kappa\Bigr)\\ 
 &= \hat{f}\,\bigl((\nabla^D_\psi s)\otimes\kappa\bigr) + \hat{f}\,(s\otimes\mathcal{L}_X\kappa)\\
 &= -\delta_{H_f\oplus~df}\bigl((\nabla^D_\psi s)\otimes\kappa\bigr) 
    - 2\pi\,\sqrt{-1}\,f\bigl((\nabla^D_\psi s)\otimes\kappa\bigr) 
    - \delta_{H_f\oplus~ df}\bigl(s\otimes \mathcal{L}_X\kappa\bigr) 
    - 2\pi\,\sqrt{-1}\,f(s\otimes \mathcal{L}_X\kappa)\\ 
 &= -(\nabla^D_{H_f\oplus~ df}\circ\nabla^D_\psi\,s)\otimes\kappa 
      - (\nabla^D_\psi s)\otimes\mathcal{L}_{H_f}\kappa 
      - 2\pi\,\sqrt{-1}\,f\bigl((\nabla^D_\psi s)\otimes\kappa\bigr)\\ 
 &\qquad - (\nabla^D_{H_f\oplus~df}s)\otimes \mathcal{L}_X\kappa 
  - s\otimes (\mathcal{L}_{H_f}\circ\mathcal{L}_X\,\kappa)
  - 2\pi\,\sqrt{-1}\,f(s\otimes \mathcal{L}_X\kappa)\\
\end{align*}
It follows from $\langle X\oplus \xi,\, H_f\oplus df\rangle_+ = 0$ that 
\[
\Lambda\bigl(X\oplus \xi,H_f\oplus df\bigr) = \frac{1}{2}\,\bigl(\,\xi(H_f) - df(X)\,) = -Xf. 
\] 
By using (\ref{sec4:formula}) and (\ref{sec3:prequantization})
these equations yield 
\begin{align*}
 (\delta_\psi\circ\hat{f} - \hat{f}\circ\delta_\psi)\,(s\otimes\kappa) 
 &= 
  -\, \bigl\{\,(\nabla^D_\psi\circ\nabla^D_{H_f\oplus~df} 
  - \nabla^D_{(H_f\oplus~df)}\circ\nabla^D_\psi)\,s\,\bigr\}\otimes\kappa\\
 &\qquad - s\otimes (\mathcal{L}_X\circ\mathcal{L}_{H_f}-\mathcal{L}_{H_f}\circ\mathcal{L}_X)\kappa 
  -  2\pi\,\sqrt{-1}\,(Xf)s\otimes\kappa\\
 &= - \left\{\Bigl(R^D_\nabla\bigl(\psi,H_f\oplus df\bigr)
    + \nabla^D_{\llbracket\psi, H_f\oplus~df\rrbracket}\Bigr)\,s\right\}\otimes \kappa\\ 
 &\qquad\qquad\qquad - s\otimes \mathcal{L}_{[X, H_f]}\kappa -  2\pi\sqrt{-1}\,(Xf)s\otimes\kappa\\ 
 &= -2\pi\,\sqrt{-1}\,\left\{\Bigl(\Lambda\bigl(X\oplus \xi,\,H_f\oplus~df\bigr) + (Xf)\Bigr)s\right\}\otimes\kappa \\
  &\qquad\qquad\qquad - \bigl(\nabla^D_{\llbracket\psi,\, H_f\oplus~df\rrbracket}\,s\bigr)\otimes \kappa 
    - s\otimes \mathcal{L}_{[X, H_f]}\kappa\\ 
 &= - \delta_{\llbracket\psi,\,H_f\oplus df\rrbracket}(s\otimes\kappa)
\end{align*}
This completes the proof.
\qquad\qquad\qquad\qquad\qquad\qquad\qquad\qquad\qquad\qquad\qquad\qquad\qquad\qquad\qquad $\Box$ 
\vspace*{0.5cm}

For a polarization $\mathcal{P}$, we define the subalgebra 
$S(\mathcal{P})$ of $(C_\mathrm{adm}^\infty(M,D),\,\{\cdot,\,\cdot\})$ as 
\[
 S({\mathcal{P}}) := \bigl\{\,f\in C^\infty_{\rm adm}(M, D)\,
 \bigm|\,\forall \psi\in\varGamma^\infty(M,\mathcal{P}):
\llbracket H_f\oplus df,\, \psi\rrbracket\in\varGamma^\infty(M, \mathcal{P})~ \,\bigr\}. 
\]
On the other hand, we define a subset $\mathfrak{H}_0$ of $\varGamma^\infty(M, L\otimes {\rm H}^{1/2})$ as 
\[
 \mathfrak{H}_0 := \{\,s\otimes\kappa\in\varGamma^\infty(M, L\otimes {\rm H}^{1/2})\,|\,
   \forall\psi\in\varGamma^\infty(M, \mathcal{P}): 
      \delta_\psi(s\otimes\kappa)=0\,\} 
\]
and assume that $\mathfrak{H}_0\ne\{\boldsymbol{0}\}$. 
Obviously, $\mathfrak{H}_0$ is a linear space over $\mathbb{C}$, and moreover 
it follows from Lemma \ref{sec5:lemma} that $\hat{f}(\mathfrak{H}_0)\subset\mathfrak{H}_0$ 
for every $f\in S(\mathcal{P})$. This leads us to the following result:
\begin{thm}
If $\mathfrak{H}_0\ne\{\boldsymbol{0}\}$, a map\, $\hat{} : S(\mathcal{P})\rightarrow 
\mathrm{End}_\mathbb{C}\,(\mathfrak{H}_0)$ by 
\[
 f\longmapsto 
 \bigl\{s\otimes \kappa \mapsto -\delta_{H_f\oplus df}(s\otimes \kappa) 
    - 2\pi\sqrt{-1}f\,(s\otimes\kappa)\bigr\}. 
\]
is a representation of $S(\mathcal{P})$ on $\mathfrak{H}_0$.
\end{thm}
We proceed with the discussion in the following two cases. 
\subsection{Compact case}
Suppose that $M$ is compact. The linear space $\mathfrak{H}_0$ has 
the inner product $\langle\cdot,\,\cdot\rangle$ defined as 
\[
 \langle s_1\otimes\kappa_1,\,s_2\otimes\kappa_2\rangle := 
 \int_M\,h(s_1,s_2)\,\overline{\kappa_1}\kappa_2 
\]
for every $s_1\otimes\kappa_1,\,s_2\otimes\kappa_2\in\mathfrak{H}_0$. 
By taking the completion of $\mathfrak{H}_0$, one obtains a Hilbert space $\mathfrak{H}$. 
The operator $\sqrt{-1}\hat{f}$ for $f\in S(\mathcal{P})$ turns out to be self-adjoint 
with respect to $\langle\cdot,\,\cdot\rangle$. Indeed, we have that 
\begin{align*}
 &\langle \hat{f}(s_1\otimes\kappa_1),\,s_2\otimes\kappa_2\rangle 
 + \langle s_1\otimes\kappa_1,\,\hat{f}(s_2\otimes\kappa_2)\rangle \\
=&\, \langle \hat{f}(s_1)\otimes\kappa_1 - s_1\otimes\mathcal{L}_{H_f}\kappa_1,
       \,s_2\otimes\kappa_2\rangle 
 + \langle s_1\otimes\kappa_1,\,
 \hat{f}(s_2)\otimes\kappa_2 - s_2\otimes\mathcal{L}_{H_f}\kappa_2\rangle\\
=&\, \langle \hat{f}(s_1)\otimes\kappa_1,\,s_2\otimes\kappa_2\rangle - 
 \langle s_1\otimes\mathcal{L}_{H_f}\kappa_1,\,s_2\otimes\kappa_2\rangle 
 + \langle s_1\otimes\kappa_1,\,\hat{f}(s_2)\otimes\kappa_2\rangle 
 - \langle s_1\otimes\kappa_1,\,s_2\otimes\mathcal{L}_{H_f}\kappa_2 \rangle \\
=& \int_M\,\Bigl\{h(\hat{f}(s_1),\,s_2)
 + h(s_1,\,\hat{f}(s_2))\Bigr\}\,\overline{\kappa_1}\kappa_2 
 - \int_M\, h(s_1,\,s_2)\,(\mathcal{L}_{H_f}\overline{\kappa_1})\,\kappa_2 
 - \int_M\, h(s_1,\,s_2)\,\overline{\kappa_1}(\mathcal{L}_{H_f}\kappa_2)\\ 
=& - \int_M\, H_f\bigl(h(s_1,\,s_2)\bigr) \,\overline{\kappa_1}\kappa_2 
- \int_M\, h(s_1,\,s_2)\,(\mathcal{L}_{H_f}\overline{\kappa_1})\,\kappa_2 
 - \int_M\, h(s_1,\,s_2)\,\overline{\kappa_1}(\mathcal{L}_{H_f}\kappa_2)\\ 
=& - \int_M\,\mathcal{L}_{H_f}\Bigl(h(s_1,\,s_2)\,\overline{\kappa_1}\kappa_2\Bigr) 
\end{align*}
for any $s_1\otimes\kappa_1,\,s_2\otimes\kappa_2\in\mathfrak{H}$. 
In the last equality, the symbol $\mathcal{L}_{H_f}(\cdot)$ denotes a Lie derivative of a 
half-density. We refer to \cite{Ythe57} for the tensor analysis of Lie derivatives. 
According to \cite{Vcoo79}, it holds that 
\[
 \int_M\,\mathcal{L}_V\kappa = 0
\]
for every vector field $V$ and every density $\kappa$ on $M$. 
From this it follows that 
\[
 \langle \hat{f}(s_1\otimes\kappa_1),\,s_2\otimes\kappa_2\rangle 
 + \langle s_1\otimes\kappa_1,\,\hat{f}(s_2\otimes\kappa_2)\rangle = 0. 
\]
This directly implies that $\sqrt{-1}\hat{f}$ is a self-adjoint operator. 
Then, the condition (\ref{sec4:bracket preserving}) holds up to the constant 
factor $\sqrt{-1}$. 
\subsection{Non-compact case}

Suppose that $M$ is not compact. We let $\mathcal{Q}$ be the subbundle of $D$ such that 
\begin{equation}\label{sec5:polarization}
 \mathcal{Q}^\mathbb{C} = \mathcal{P}\cap\overline{\mathcal{P}}
\end{equation}
and assume that a singular distribution $M\ni p\mapsto (\rho_{TM})_p(\mathcal{Q}_p)\subset T_pM$ 
defines a regular foliation $\mathcal{F}$ whose leaf space $N=M/\mathcal{F}$ is a Hausdorff manifold. 
We denote by $q_{\mathcal{F}}$ the natural projection from $M$ to $N$. 
For any $f\in S(\mathcal{P})$ and $X\oplus \xi\in\mathcal{Q}$, 
the vector field $\llbracket\,H_f\oplus df,\,X\oplus \xi\,\rrbracket$ is tangent to each $q_\mathcal{F}$-fiber: 
$\llbracket\,H_f\oplus df,\,X\oplus \xi\,\rrbracket = [H_f,X]\oplus \mathcal{L}_{H_f}\xi
\in\varGamma^\infty(M,\mathcal{Q})$. Accordingly, it holds that 
$[H_f,X]\in \rho_{TM}\,(\mathcal{Q})$ for any vector field $X\in\rho_{TM}\,(\mathcal{Q})$. 
This means that $H_f$ is a lift of a smooth vector field on $N$. 
Let ${\rm H}_N^{1/2}$ be a half density bundle over $N$ and 
$\kappa_N$ a half density. Then, from (\ref{sec5:representation}) it holds that 
\[
 \hat{f}(s\otimes(q_\mathcal{F})^*\kappa_N) = (\hat{f}s)\otimes(q_\mathcal{F})^*\kappa_N 
- s\otimes(\mathcal{L}_{H_f}(q_\mathcal{F})^*\kappa_N)
=(\hat{f}s)\otimes(q_\mathcal{F})^*\kappa_N - s\otimes(q_\mathcal{F})^*(\mathcal{L}_{{q_\mathcal{F}}_*(H_f)}\kappa_N) 
\]
for any $s\otimes(q_\mathcal{F})^*\kappa_N\in\varGamma^\infty(M,L\otimes (q_\mathcal{F})^*{\rm H}_N^{1/2})$. 
This enables us to consider the half densities of $M$ which are transversal to $\mathcal{F}$. 
Since $\mathcal{L}_{X}((q_\mathcal{F})^*\kappa_N)=(q_\mathcal{F})^*(\mathcal{L}_{\pi_*X}\kappa_N)=0$ 
for $(X,\xi)\in\mathcal{Q}$, we have that 
\begin{equation*}
 \mathcal{L}_{X}\Bigl(\bigl((q_\mathcal{F})^*\kappa^1_N\bigr)\otimes\bigl((q_\mathcal{F})^*\kappa^2_N\bigr)\Bigr)
 = (q_\mathcal{F})^*(\mathcal{L}_{{q_\mathcal{F}}_*X}\kappa^1_N)\otimes(q_\mathcal{F})^*\kappa^2_N 
  + (q_\mathcal{F})^*\kappa^1_N\otimes(q_\mathcal{F})^*(\mathcal{L}_{{q_\mathcal{F}}_*X}\kappa^2_N)
 = 0 
\end{equation*}
for every $(q_\mathcal{F})^*\kappa^1_N, (q_\mathcal{F})^*\kappa^2_N\in (q_\mathcal{F})^*{\rm H}^{1/2}_N$. 
In other words, the tensor field $((q_\mathcal{F})^*\kappa^1_N)\otimes((q_\mathcal{F})^*\kappa^2_N)$ on $M$ is 
invariant under the flow of $X\in\rho_{TM}(\mathcal{Q})$. 
Therefore, there exists a 1-density $\nu_N$ of $N$ onto which 
$((q_\mathcal{F})^*\kappa^1_N)\otimes((q_\mathcal{F})^*\kappa^2_N)$ projects. 
As a result, we let $\mathfrak{H}_1$ be the linear subspace of $\mathfrak{H}_0$ consisting of 
the elements in $\varGamma^\infty(M,L\otimes (q_\mathcal{F})^*{\rm H}_N^{1/2})$ which have compact support 
in $N$, and assume that $\mathfrak{H}_1\ne \{\boldsymbol{0}\}$. 
We define the inner product $\langle\cdot,\cdot\rangle$ on $\mathfrak{H}_1$ as 
\[
 \langle s_1\otimes(q_\mathcal{F})^*\kappa_N^1,\,s_2\otimes(q_\mathcal{F})^*\kappa_N^2\rangle := 
 \int_N\,h(s_1,s_2)\,\nu_N 
\]
for every $s_1\otimes(q_\mathcal{F})^*\kappa^1_N,\,s_2\otimes (q_\mathcal{F})^*\kappa^2_N\in\mathfrak{H}_1$. 
Replacing $\mathfrak{H}_0$ with $\mathfrak{H}_1$, we obtain a Hilbert space from $\mathfrak{H}_1$ 
and find that $\sqrt{-1}\hat{f}$ for $f\in S(\mathcal{P})$ is a self-adjoint operator on 
$\mathrm{End}_{\mathbb{C}}\,(\mathfrak{H}_1)$ 
in a way similar to the compact case. 
\vspace*{0.5cm}

\noindent{\bf Acknowledgments.}~The author would like to express his deepest 
gratitude to Emeritus Professor Toshiaki Kori of Waseda University for 
useful discussion and various support. 
He also would appreciate Professor Alan Weinstein and the referees very much for helpful comments, 
pointing out typos and suggesting many improvements to the first version of the manuscript. 
Lastly, he wishes to thank Keio University and Waseda University for the hospitality while part of the work was being done.

\noindent\footnotesize{
{\em Yuji Hirota}\\
\texttt{School of Veterinary Medicine, Azabu University, JAPAN 
\\ hirota@azabu-u.ac.jp}}
\end{document}